    \def\version{January 6, 2003}

\newif\ifpdf
\ifx\pdfoutput\undefined
\pdffalse 
\else
\pdfoutput=1 
\pdftrue \fi

\newif\ifdate
\datetrue
\datefalse

\documentclass[reqno,twoside,11pt]{amsart}
\usepackage{cite}
\usepackage{amsmath}
\usepackage{amsfonts}
\usepackage{amssymb}
\usepackage{epsf}
\usepackage{graphicx}
\usepackage{graphics}
\usepackage{verbatim}

\IfFileExists{myowntimes.sty}
	{\usepackage{myowntimes}}{\usepackage{times}\usepackage{mathrsfs}}

\DeclareFontFamily{OT1}{eusm}{} \DeclareFontShape{OT1}{eusm}{m}{n}
{<5> <6> <7> <8> <9> <10> <11> <12> <14.4> eusm10}{}
\DeclareMathAlphabet{\eusm}{OT1}{eusm}{m}{n}

\newtheoremstyle{thm}{1.5ex}{1.5ex}{\itshape\rmfamily}{}
{\bfseries\rmfamily}{}{2ex}{}

\newtheoremstyle{rem}{1.3ex}{1.3ex}{\rmfamily}{}
{\itshape}
{} {1.5ex}{}

\newenvironment{proofsect}[1]
{\vskip0.1cm\noindent{\rmfamily\itshape
#1.}}{\qed\vspace{0.15cm}}

\theoremstyle{thm}
\newtheorem{theorem}{Theorem}[section]
\newtheorem{lemma}[theorem]{Lemma}
\newtheorem{proposition}[theorem]{Proposition}

\newtheorem*{Main Theorem}{Main Theorem.}
\newtheorem*{Key Estimate}{Key Estimate}
\newtheorem{corollary}[theorem]{Corollary}

\newtheorem{property}{Property}
\newtheorem*{Uniformity Property}{Uniformity Property}

\theoremstyle{rem}

\numberwithin{equation}{section}

\renewcommand{\section}{\secdef\sct\sect}
\newcommand{\sct}[2][default]{\refstepcounter{section}
\addcontentsline{toc}{section}
{{\tocsection {}{\thesection}{\!\!\!\!#1\dotfill}}{}}
\vspace{0.7cm}
\centerline{ 
\scshape\arabic{section}.\ #1} \nopagebreak \vspace{0.2cm}}
\newcommand{\sect}[1]{
\vspace{0.4cm} \centerline{\large\scshape\rmfamily #1}
\vspace{0.2cm}}

\renewcommand{\subsection}{\secdef\subsct\sbsect}
\newcommand{\subsct}[2][default]{\refstepcounter{subsection}
\addcontentsline{toc}{subsection}
{{\tocsection{\!\!}{\hspace{1.2em}\thesubsection}{\!\!\!\!#1\dotfill}}{}}
\nopagebreak \vspace{0.45\baselineskip} {\flushleft\bf
\arabic{section}.\arabic{subsection}~\bf #1.~}
\\*[3mm]\noindent
\nopagebreak}
\newcommand{\sbsect}[1]{\vspace{0.1cm}\noindent
\textbf{#1.~}\vspace{0.1cm}}

\renewcommand{\subsubsection}{%
\secdef \subsubsect\sbsbsect}
\newcommand{\subsubsect}[2][default]{%
\refstepcounter{subsubsection} 
\addcontentsline{toc}{subsubsection}{{\tocsection{\!\!}
{\hspace{3.05em}\thesubsubsection}{\!\!\!\!#1\dotfill}}{}}
\nopagebreak
\vspace{0.15\baselineskip} \nopagebreak {\flushleft\rmfamily
\itshape\arabic{section}.\arabic{subsection}.\arabic{subsubsection}
\ \rmfamily #1\/.}\ }
\newcommand{\sbsbsect}[1]{\vspace{0.1cm}\noindent
\rmfamily \itshape
\arabic{section}.\arabic{subsection}.\arabic{subsubsection} \
\sffamily #1\/.\ }

\newcommand{\dist}{\operatorname{dist}}

\newcommand{\textd}{\text{\rm d}}

\newcommand{\diag}{\text{\rm diag}}

\newcommand{\aff}{{\text{\rm aff}\,}}
\newcommand{\ri}{{\text{\rm ri}\,}}

\newcommand{\GG}{\mathcal G}

\newcommand{\UU}{\mathcal U}

\newcommand{\E}{\mathbb E}

\newcommand{\K}{\mathbb K}

\newcommand{\R}{\mathbb R}
\newcommand{\BbbS}{\mathbb S}

\newcommand{\V}{\mathbb V}

\newcommand{\Z}{\mathbb Z}

\newcommand{\1}{{\text{\sf 1}}}

\newcommand{\bS}{{\boldsymbol S}}
\newcommand{\bQ}{{\boldsymbol Q}}

\newcommand{\bm}{{\boldsymbol m}}
\newcommand{\bk}{{\boldsymbol k}}
\newcommand{\bM}{{\boldsymbol M}}
\newcommand{\ba}{{\boldsymbol a}}
\newcommand{\bb}{{\boldsymbol b}}
\newcommand{\bh}{{\boldsymbol h}}
\newcommand{\bg}{{\boldsymbol g}}

\newcommand{\bs}{{\boldsymbol s}}
\newcommand{\bv}{{\boldsymbol v}}
\newcommand{\bw}{{\boldsymbol w}}
\newcommand{\bzero}{{\boldsymbol 0}}
\newcommand{\blambda}{{\boldsymbol\lambda}}
\newcommand{\bomega}{{\boldsymbol\omega}}
\newcommand{\btau}{{\boldsymbol\tau}}

\newcommand{\MF}{{\text{\rm MF}}}
\newcommand{\conv}{{\text{\rm Conv}}}

\newcommand{\hate}{{\text{\rm \^e}}}
\newcommand{\hatv}{{\text{\rm \^v}}}

\newcommand{\eusmJ}{{\eusm J}}

\newcommand{\twoeqref}[2]{(\ref{#1}--\ref{#2})}
\newcommand{\verzicka}{\ifdate, \version\fi}

\begin{document}

\title[Phase transitions and mean-field theory\verzicka]
{
\fontsize{15}{20}\selectfont 
Rigorous analysis of discontinuous phase~transitions via
mean-field bounds}
\author[Marek~Biskup and Lincoln~Chayes\verzicka]
{Marek~Biskup\, and\, Lincoln~Chayes}
\maketitle

\thispagestyle{empty}

\vspace{-5mm} 
\centerline{\textit{Department of Mathematics, UCLA, Los Angeles CA 90095-1555, U.S.A.}}

\vspace{1mm}
\begin{quote}
{\footnotesize \textbf{Abstract:} } \footnotesize
We consider a variety of nearest-neighbor spin models defined on the $d$-dimensional hypercubic lattice~$\Z^d$. Our essential assumption is that these models satisfy the condition of reflection positivity. We prove that whenever the associated mean-field theory predicts a discontinuous transition, the actual model also undergoes a discontinuous transition (which occurs near the mean-field transition temperature), provided the dimension is sufficiently large or the first-order transition in the mean-field model is sufficiently strong. As an application of our general theory, we show that for~$d$ sufficiently large, the~$3$-state Potts ferromagnet on~$\Z^d$ undergoes a first-order phase transition as the temperature varies.  Similar results are established for all~$q$-state Potts models with $q\ge3$, the~$r$-component cubic models with $r\ge4$ and the~$O(N)$-nematic liquid-crystal models with~$N\ge3$.
\end{quote}



\vspace{-0.3cm}

\setcounter{tocdepth}{3}
\footnotesize
\begin{list}{}
{\setlength{\topsep}{0in}\setlength{\leftmargin}{0.34in}\setlength{\rightmargin}{0.5in}}
\item[]
\tableofcontents
\end{list}
\vspace{-1.3cm}
\begin{list}{}
{\setlength{\topsep}{0in}\setlength{\leftmargin}{0.5in}\setlength{\rightmargin}{0.5in}}
\item[]
\hskip-0.01in
\hbox to 2.3cm{\hrulefill}
\item[]
{\fontsize{8.5}{8.5}\selectfont\copyright\,\,\,Copyright rests with the authors. Reproduction
of the entire article for non-commercial purposes
is permitted without charge.\vspace{2mm}}
\end{list}
\normalsize

\section{Introduction}
\label{sec1}\vspace{-0.2cm}
\subsection{Motivation and outline}
\label{sec1.1}\noindent
Mean-field theory has traditionally played a seminal role for qualitative understanding of phase transitions. In fact, most practical studies of complex physical systems begin (and sometimes end) with the analysis of the corresponding mean field theory. The central idea of mean-field theory---dating back to \cite{Curie,Weiss}---is rather compelling: The ostensibly complicated interactions acting on a particular element of the system are replaced by the action of an effective (or \emph{mean}) external field. This field causes a response at the point of question and its value has to be self-consistently adjusted so that the response matches the effective field. The practical outcome of this procedure is a set of equations, known as the \emph{mean-field equations}. In contrast to the original, fully interacting system, the mean-field equations are susceptible to direct analytical or~numerical~methods.

There is a general consensus that mean-field predictions are qualitatively or even quantitatively accurate. However, for short-range systems, a mathematical foundation of this belief has not been presented in a general context. A number of rigorous results have related various lattice systems to their mean-field counterparts, either in the form of bounds on transition temperatures and critical exponents, see~\cite{Ellis,Simon,FFS} and references therein, or in terms of limits of the free energy~\cite{PT} and the magnetization~\cite{BKLS,KS} as the dimension tends to infinity. In all of these results, the nature of the phase transition is not addressed or the proofs require special symmetries which, as it turns out, ensure that the transition is continuous. But, without special symmetries (or fine tuning) phase transitions are typically discontinuous, so generic short-range systems have heretofore proved elusive. (By contrast, substantial progress along these lines has been made for systems where the range of the interaction plays the role of a large parameter. See, e.g., \cite{Cassandro-Presutti,Lebowitz-Mazel-Presutti,Bovier-Zahradnik1,
Bovier-Zahradnik2}.)

In this paper we demonstrate that for a certain class of nearest-neighbor spin systems, namely those that are \emph{reflection positive}, mean-field theory indeed provides a rigorous guideline for the order of the transition. In particular, we show that the actual systems undergo a first-order transition whenever the associated mean-field model predicts this behavior, provided the spatial dimension is sufficiently high and/or the phase transition is sufficiently strong. Furthermore, we give estimates on the difference between the values of parameters of the actual model and its mean-field counterpart at their corresponding transitions and show that these differences tend to zero as the spatial dimension tends to infinity. In short, mean field theory is \emph{quantitatively accurate} whenever the dimension is sufficiently large. 

The main driving force of our proofs is the availability of the so called \emph{infrared bound}~\cite{FSS,DLS,FILS1,FILS2}, which we use for estimating the correlations between nearest-neighbor spins. It is worth mentioning that the infrared bound is the principal focus of interest in a class of rigorous results on mean-field \emph{critical} behavior of various combinatorial models~\cite{Brydges-Spencer,Hara-Slade1,Hara-Slade2,Hara-Slade3,Remco-Frank-Gordon,Remco-Gordon} and percolation~\cite{Hara-Slade4,Hara-Slade5,Hara-Slade6,Hara-Slade7,Remco-Frank-Gordon2,Remco-Slade1a,Hara-Remco-Slade2} based on the technique of the lace expansion. However, in contrast to these results (and to the hard work that they require), our approach is more reminiscent of the earlier works on high-dimensional systems~\cite{Aizenman,AF,ABF}, where the infrared bound is provided as an \emph{input}. In particular, for our systems this input is a consequence of reflection positivity. (As such, some of our results can also be extended to systems with long-range forces; the relevant modifications will appear in a separate publication~\cite{Biskup-Chayes}.)

The principal substance of this paper is organized as follows: We devote the remainder of Section~\ref{sec1} to a precise formulation of the general class of spin systems that we consider, we then develop some general mean-field formalism and, finally, state our main theorems. Section~\ref{sec2} contains a discussion of three eminent models---Potts, cubic and nematic---with specific statements of theorems which underscore the first-order (and mean-field) nature of the  phase transitions for the \hbox{large-$d$} version of these models. In Section~\ref{sec3} we develop and  utilize  the principal tools needed in this work and provide proofs of all statements made in Section~\ref{sec1}. In Section~\ref{sec4}, we perform detailed analyses and collect various known results  on the mean-field theories for the specific models mentioned above. When these systems are ``sufficiently prepared,'' we apply the Main Theorem to prove all of the results stated in Section~\ref{sec2}. Finally, in Section~\ref{sec5}, we show that for any model in the class considered, the mean-field theory can be realized by defining the problem on the complete graph.

\subsection{Models of interest}
\label{sec1.2}\noindent
Throughout this paper, we will consider the following class of spin systems on the~$d$-dimensional hypercubic lattice~$\Z^d$: The \emph{spins}, denoted by $\bS_x$, take values in some fixed set~$\Omega$, which is a subset of a finite dimensional vector space~$\E_\Omega$. 
We will use $(\cdot\,,\cdot)$ to denote the (positive-definite) inner product in $\E_\Omega$ and assume that $\Omega$ is compact in the topology induced by this inner product. The spins are weighted according to an \emph{a priori} Borel probability measure~$\mu$ whose support is~$\Omega$. An assignment of a spin value~$\bS_x$ to each site~$x\in\Z^d$ defines a \emph{spin configuration}; we assume that the \emph{a priori} joint distribution of all spins on~$\Z^d$ is i.i.d. Abusing the notation slightly, we will use~$\mu$ to denote the joint \emph{a priori} measure on spin configurations and use~$\langle-\rangle_0$ to denote the expectation with respect to~$\mu$. 

The interaction between the spins is described by the (formal) Hamiltonian
\begin{equation}
\label{Ham} \beta H = -\frac J{2d}\sum_{\langle
x,y\rangle}(\bS_x,\bS_y)-\sum_x(\bb,\bS_x).
\end{equation}
Here $\langle x,y\rangle$ denotes a nearest-neighbor pair of~$\Z^d$, the quantity~$\bb$, playing the role of an external field, is a vector from~$\E_\Omega$ and~$\beta$, the inverse temperature, has been incorporated into the (normalized) coupling constant~$J\ge0$ and the field parameter~$\bb$. 

The interaction Hamiltonian gives rise to the concept of a Gibbs measure which is defined as follows: Given a finite set $\Lambda\subset\Z^d$, a configuration $\bS=(\bS_x)_{x\in\Lambda}$ in $\Lambda$ and a boundary condition $\bS'=(\bS'_x)_{x\in\Z^d\setminus\Lambda}$ in $\Z^d\setminus\Lambda$, we let $\beta H_\Lambda(\bS|\bS')$ be given by \eqref{Ham} with the first sum on the right-hand side of \eqref{Ham} restricted to $\langle x,y\rangle$ such that $\{x,y\}\cap\Lambda\ne\emptyset$, the second sum restricted to~$x\in\Lambda$, and~$\bS_x$ for~$x\not\in\Lambda$ replaced by~$\bS_x'$. Then we define  the measure $\nu_\Lambda^{(\bS')}$ on configurations~$\bS$ in~$\Lambda$ by the expression
\begin{equation}
\label{nuSS}
\nu_\Lambda^{(\bS')}(\textd\bS)=\frac{e^{-\beta H_\Lambda(\bS|\bS')}}{Z_\Lambda(\bS')}\mu(\textd\bS),
\end{equation}
where~$Z_\Lambda(\bS')$ is the appropriate normalization constant which is called the \emph{partition function}. The measure in \eqref{nuSS} is the \emph{finite-volume Gibbs measure} corresponding to the interaction \eqref{Ham}.

In statistical mechanics, the measure \eqref{nuSS} describes the thermodynamic equilibrium of the spin system in~$\Lambda$. To address the question of phase transitions, we have to study the possible limits of these measures as~$\Lambda$ expands to fill in~$\Z^d$. In accord with the standard definitions, see~\cite{Georgii}, we say that the spin model undergoes a \emph{first-order phase transition} at parameter values $(J,\bb)$ if there are at least two distinct infinite-volume limits of the measure in \eqref{nuSS} arising from different boundary conditions. We will call these limiting objects either infinite-volume Gibbs measures or, in accordance with mathematical-physics nomenclature, \emph{Gibbs states}. We refer the reader to \cite{Georgii,Simon} for more details on the general properties of Gibbs states and phase transitions.

We remark that, while the entire class of models has been written so as to appear identical, the physics will be quite different depending on the particulars of~$\Omega$ and~$\mu$, and the inner product. Indeed, the language of magnetic systems has been adapted only for linguistic and notational convenience. The above framework can easily accommodate any number of other physically motivated interacting models such as lattice gases, ferroelectrics, etc.

\subsection{Mean-field formalism}
\label{sec1.3}\noindent
Here we will develop the general formalism needed for stating the principal mean-field bounds. The first object of interest is the logarithmic moment generating function of the distribution~$\mu$,
\begin{equation}
\label{G} 
G(\bh)=\log\int_\Omega \mu(\textd\bS)\, e^{(\bS,\bh)}.
\end{equation}
Since~$\Omega$ was assumed compact, $G(\bh)$ is finite for all $\bh\in\E_\Omega$. Moreover, $\bh\mapsto G(\bh)$ is continuous and convex throughout~$\E_\Omega$.

Every mean-field theory relies on a finite number of thermodynamic functions of internal responses. For the systems with interaction~\eqref{Ham}, the object of principal interest is the \emph{magnetization}. In general, magnetization is a quantity taking values in the closed, convex hull of $\Omega$, here denoted by $\conv(\Omega)$. If $\bm\in\conv(\Omega)$, then the \emph{mean-field entropy function} is defined via a Legendre transform of~$G(\bh)$,
\begin{equation}
\label{S} S(\bm)=\inf_{\bh\in\E_\Omega}\bigl\{G(\bh)-(\bm,\bh)\bigr\}.
\end{equation}
(Strictly speaking, \eqref{S} makes sense even for $\bm\not\in\conv(\Omega)$ for which we simply get $S(\bm)=-\infty$.) In general, $\bm\mapsto S(\bm)$ is concave and we have $S(\bm)\le0$ for all $\bm\in\conv(\Omega)$. From the perspective of the large-deviation theory (see~\cite{Ellis,Dembo-Zeitouni}), the mean-field entropy function is (the negative of) the rate function for the probability that the average of many spins is near~$\bm$.

To characterize the effect of the interaction, we have to introduce energy into the game.
For the quadratic Hamiltonian in \eqref{Ham}, the \emph{(mean-field) energy function} is given simply by
\begin{equation}
\label{E} E_{J,\bb}(\bm)=-\frac12 J|\bm|^2-(\bm,\bb),
\end{equation}
where $|\bm|^2=(\bm,\bm)$. On the basis of physical considerations, a state of thermodynamic equilibrium corresponds to a balance between the energy and the entropy. The appropriate thermodynamic function characterizing this balance is the free energy. We therefore define the \emph{mean-field free-energy function} by setting~$\varPhi_{J,\bb}(\bm)=E_{J,\bb}(\bm)-S(\bm)$, i.e.,
\begin{equation}
\label{Fi} \varPhi_{J,\bb}(\bm)=-\frac 12
J|\bm|^2-(\bm,\bb)-S(\bm).
\end{equation}
The mean-field (Gibbs) free energy $F_\MF(J,\bb)$ is defined by minimizing~$\varPhi_{J,\bb}(\bm)$ over all $\bm\in\conv(\Omega)$. Assuming a unique minimizer, this and (\ref{S}-\ref{E}) give us a definition of the mean-field magnetization, entropy and energy. A more interesting situation occurs when there is more than one minimizer of~$\varPhi_{J,\bb}$. The latter cases are identified as the points of phase coexistence while the former situation is identified as the uniqueness region.

For the sake of completeness, it is interesting to observe that every minimizer of~$\varPhi_{J,\bb}(\bm)$ (in fact, every stationary point) in the relative interior of $\conv(\Omega)$ is a solution of the equation
\begin{equation}
\label{MFeq}
\bm=\nabla G(J\bm+\bb),
\end{equation}
where~$\nabla$ denotes the (canonical) gradient in~$\E_\Omega$. This is the \emph{mean-field equation} for the magnetization, which describes the self-consistency constraint that we alluded to in Section~\ref{sec1.1}. The relation between \eqref{MFeq} and the stationarity of~$\varPhi_{J,\bb}$ is seen as follows: $\nabla\varPhi_{J,\bb}(\bm)=0$ implies that $J\bm+\bb+\nabla S(\bm)=0$. But $\bh=-\nabla S(\bm)$ is equivalent to $\bm=\nabla G(\bh)$, and stationarity therefore~implies~\eqref{MFeq}.

We conclude with a claim that an immediate connection of the above formalism to \emph{some} statistical mechanics problem is possible. Indeed, if the Hamiltonian \eqref{Ham} is redefined for the \text{complete graph} on~$N$ vertices, then the quantity~$\varPhi_{J,\bb}(\bm)$ emerges as the rate function in a large-deviation principle for magnetization and hence $F_\MF(J,\bb)$ \emph{is} the free energy in this model. A precise statement and a proof will appear in the last section (Theorem~\ref{thm5.1} in Section~\ref{sec5}); special cases of this result have been known since time immemorable, see e.g.~\cite{Ellis}.

\subsection{Main results}
\label{sec1.4}\noindent
Now we are in a position to state our general results. The basic idea is simply to watch what happens when the value of the magnetization in an actual system (governed by~\eqref{Ham}) is inserted into the associated mean-field free-energy function. We begin with a general bound which relies only on convexity:

\begin{theorem}
\label{thm1}
Consider the spin system on $\Z^d$ with the Hamiltonian~\eqref{Ham} and let~$\nu_{J,\bb}$ be an infinite-volume Gibbs measure corresponding to the parameters $J\ge0$ and $\bb\in\E_\Omega$ in \eqref{Ham}. Suppose that~$\nu_{J,\bb}$ is invariant under the group of translations and rotations of~$\Z^d$. Let $\langle-\rangle_{J,\bb}$ denote the expectation with respect to~$\nu_{J,\bb}$ and let~$\bm_\star$ be the magnetization of the state~$\nu_{J,\bb}$ defined by
\begin{equation}
\label{bmstar}
\bm_\star=\langle\bS_0\rangle_{J,\bb},
\end{equation}
where~$0$ denotes the origin in~$\Z^d$.
Then
\begin{equation}
\label{genbd}
\varPhi_{J,\bb}(\bm_\star)\le
\inf_{\bm\in\conv(\Omega)}\varPhi_{J,\bb}(\bm)
+\frac J2
\bigl[\bigl\langle(\bS_0,\bS_x)\bigr\rangle_{J,\bb}
-|\bm_\star|^2\bigr],
\end{equation}
where~$x$ denotes a nearest neighbor of the origin.
\end{theorem}

Thus, whenever the fluctuations of nearest-neighbor spins have small correlations, the physical magnetization \emph{almost} minimizes the mean-field free energy. The bound \eqref{genbd} immediately leads to the following observation, which, to the best of our knowledge, does not  appear in the literature:

\begin{corollary}
\label{cor-ener1}
Let $\nu_{J,\bb}$ and $\langle-\rangle_{J,\bb}$ be as in Theorem~\ref{thm1} and let $\bm_\star$ be as in \eqref{bmstar}. Then
\begin{equation}
\bigl\langle(\bS_x,\bS_y)\bigr\rangle_{J,\bb}\ge
|\bm_\star|^2
\end{equation}
for any pair of nearest-neighbors~$x,y\in\Z^d$. In particular, for any model with interaction \eqref{Ham}, the nearest-neighbor spins are positively correlated in any Gibbs state which is  invariant under the translations and rotations of~$\Z^d$.
\end{corollary}

Our next goal is to characterize a class of Gibbs states for which the correlation term on the right-hand side of \eqref{genbd} is demonstrably small. However, our proofs will make some minimal demands on the Gibbs states themselves and it is therefore conceivable that we may not be able to access \emph{all} the extremal magnetizations. To define those values of magnetization for which our proofs hold, let~$F(J,\bb)$ denote the infinite-volume free energy per site of the system on $\Z^d$, defined by taking the thermodynamic limit of $-\frac1{|\Lambda|}\log Z_\Lambda$, see e.g.~\cite{Ruelle}. (Note that the existence of this limit follows automatically by the compactness of~$\Omega$.) The function~$F(J,\bb)$ is concave and, therefore, has all directional derivatives. Let $\mathscr{K}_\star(J,\bb)$ be the set of all pairs $[e_\star,\bm_\star]$ such that
\begin{equation}
\label{physF}
F(J+\Delta J,\bb+\Delta\bb)-
F(J,\bb)\le e_\star\Delta J+(\bm_\star,\Delta\bb)
\end{equation}
holds for all numbers~$\Delta J$ and all vectors~$\Delta\bb\in\E_\Omega$.
By a well-known result (see the
discussion of the properties of \emph{subdifferential} on
page~215 of~\cite{Rocky}), $\mathscr{K}_\star(J,\bb)$ is a convex set; we
let $\mathscr{M}_\star(J,\bb)$ denote the set of all values
$\bm_\star$ such that $[e_\star,\bm_\star]$ is an extreme point of
the set $\mathscr{K}_\star(J,\bb)$ for some value $e_\star$. 

\smallskip
Our Main Theorem is then as follows:

\begin{Main Theorem}
Let $d\ge 3$ and consider the spin system on $\Z^d$ with the Hamiltonian~\eqref{Ham}. Let~$n$ denote the dimension of\/ $\E_\Omega$. For $J\ge0$ and $\bb\in\E_\Omega$, let $\bm_\star\in\mathscr{M}_\star(J,\bb)$. Then
\begin{equation}
\label{5} 
\varPhi_{J,\bb}(\bm_\star)\le
\inf_{\bm\in\conv(\Omega)}\varPhi_{J,\bb}(\bm)+Jn\frac\kappa2 I_d,
\end{equation}
where $\kappa=\max_{\bS\in\Omega}(\bS,\bS)$ and
\begin{equation}
\label{Id} 
I_d=\int_{[-\pi,\pi]^d} 
\frac{\textd^d\!k}{(2\pi)^d}
\frac{[1-\widehat D(k)]^2}{\widehat D(k)}
\end{equation}
with $\widehat D(k)=1-\frac1d\sum_{j=1}^d\cos(k_y)$.
\end{Main Theorem}

\begin{figure}[t]
\refstepcounter{obrazek}
\ifpdf \centerline{\includegraphics[width=3.5truein]{Potts.pdf}}
\else
\centerline{\epsfxsize=3.5truein \epsffile{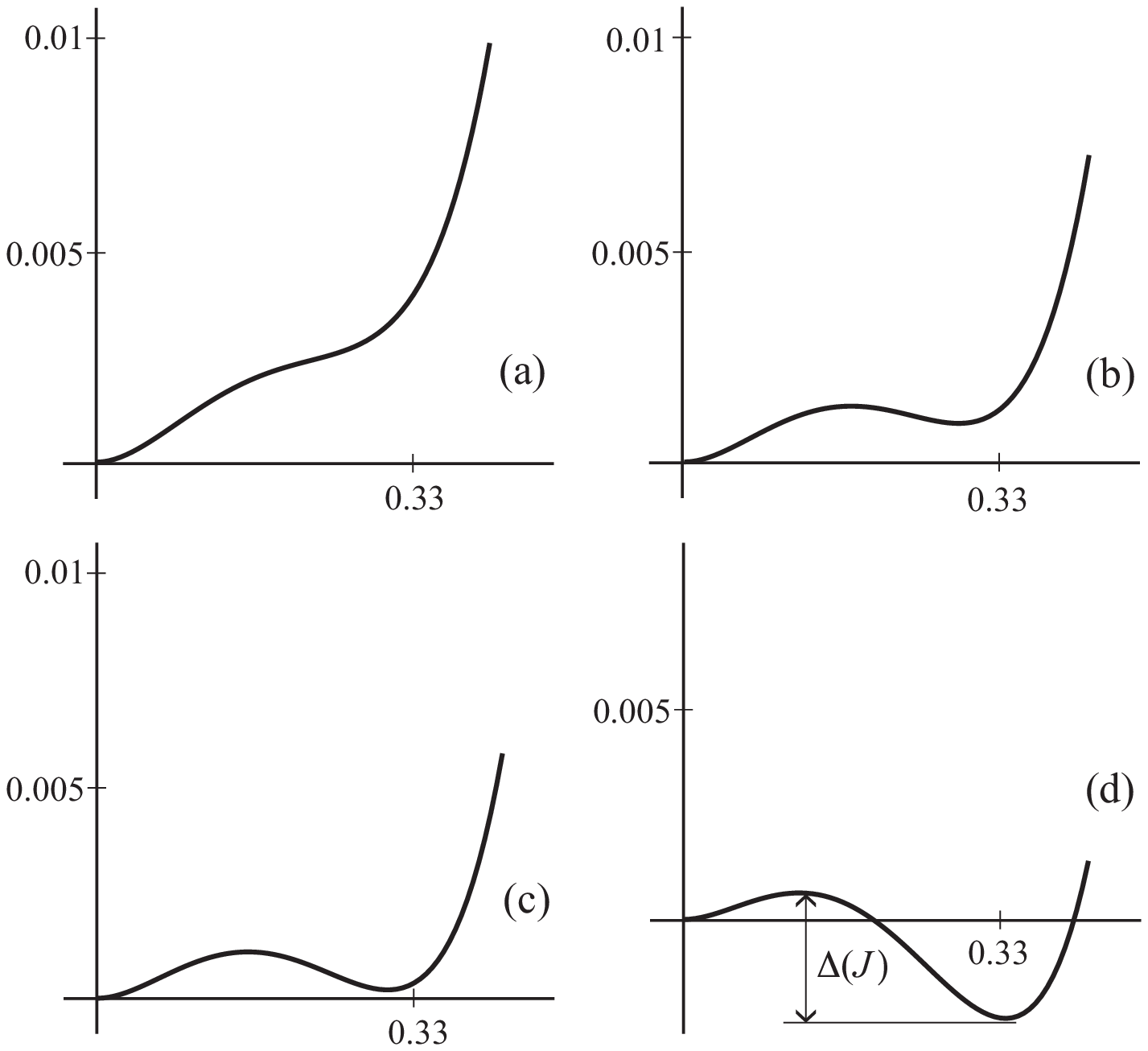}} 
\fi
\bigskip
\begin{quote}
\fontsize{9.5}{7}\selectfont
{\sc Figure~1.\ }
\label{fig1}
The mean-field free energy as a function of a scalar magnetization $m(J)$ for the typical model undergoing a first-order phase transition. In an interval of values of~$J$, there are two local minima which switch their order at $J=J_\MF$. If the ``barrier'' height $\Delta(J)$ always exceeds the error term from \eqref{5}, there is a forbidden interval of scalar magnetizations and $m(J)$ has to jump as~$J$ varies. The actual plot corresponds to the $3$-state Potts model for~$J$ taking the values (a)~$2.73$, (b)~$2.76$, (c)~$2.77$ and (d)~$2.8$. See Section~\ref{sec2.1} for more details.
\normalsize
\end{quote}
\end{figure}

\begin{figure}[t]
\refstepcounter{obrazek}
\ifpdf \centerline{\includegraphics[width=2.8truein]{magnetization.pdf}}
\else
\centerline{\epsfxsize=2.8truein \epsffile{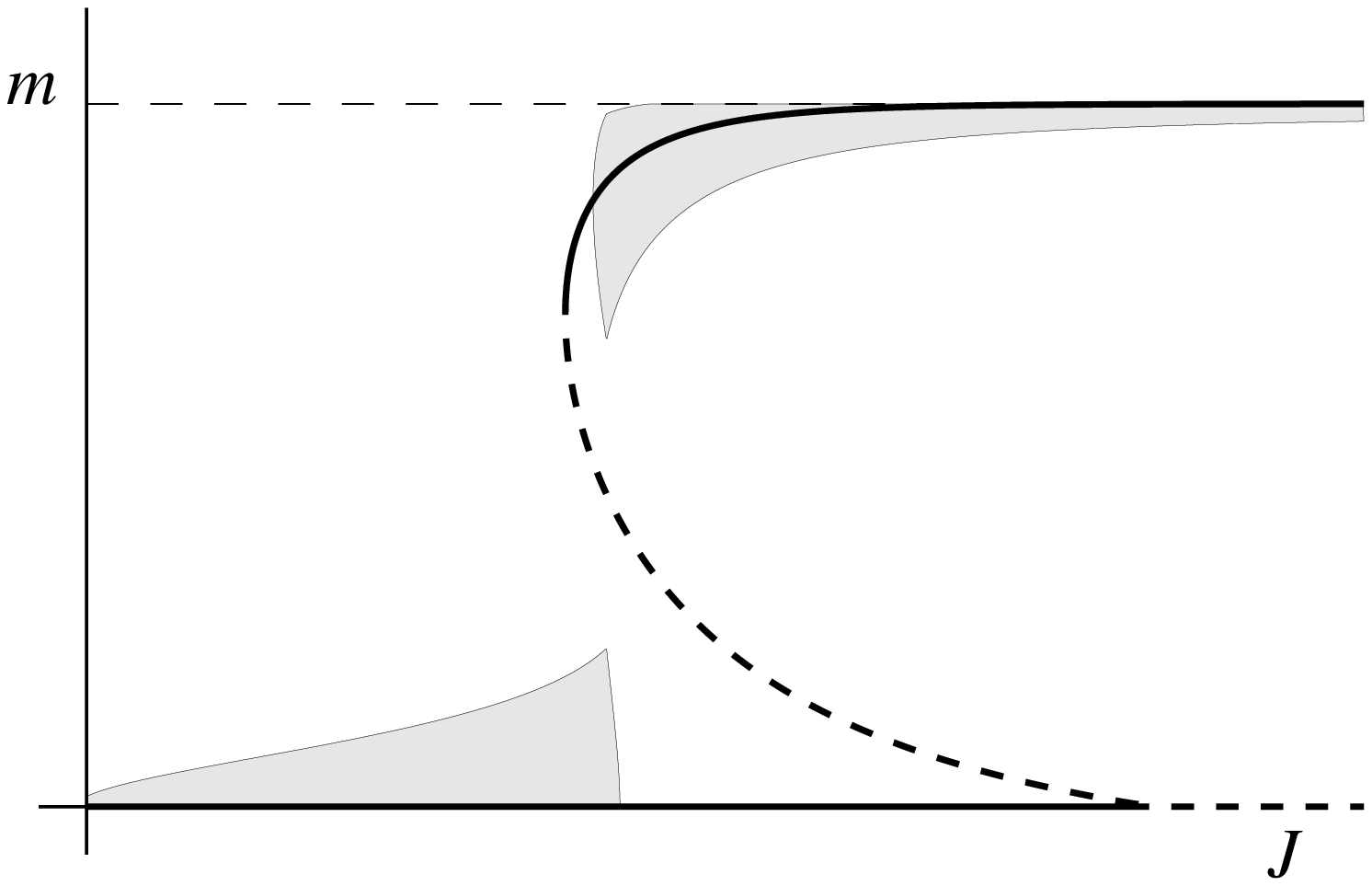}} 
\fi
\bigskip
\begin{quote}
\fontsize{9.5}{7}\selectfont
{\sc Figure~2.\ }
\label{fig2}
The solutions of the mean-field equation for the scalar order parameter~$m$ as a function of~$J$ for the $10$-state Potts model. The solid lines indicate the local minima, the dashed lines show the other solutions to the mean-field equation. The portions of these curves in the regions where $m$ is sufficiently close to zero or one can be (rigorously) controlled using perturbative calculations. These alone prove that the mean-field  theory ``does not admit continuous solutions'' and, therefore, establish a  first order transitions for~$d\gg1$. The shaded regions show the set of allowed magnetizations for the system on~$\Z^d$ when $I_d\le0.002$. In addition to manifestly proving a discontinuous transition, these provide tight numerical bounds on the transition temperature and reasonable bounds on the size of the jump.
\normalsize
\end{quote}
\end{figure}

The bound \eqref{5} provides us with a powerful method for proving first-order phase transitions on the basis of a comparison with the associated mean-field theory. The key to our whole program is that the ``error term'', $Jn\frac\kappa2 I_d$, vanishes in the $d\to\infty$ limit; in fact, 
\begin{equation}
I_d=\frac1{2d}\bigl(1+o(1)\bigr)\quad\text{as}\quad d\to\infty,
\end{equation}
see~\cite{BKLS}. For~$d$ sufficiently large, the bound~\eqref{5} thus forces the magnetization of the actual system to be \emph{near} a value of~$\bm$ that \emph{nearly} minimizes~$\varPhi_{J,\bb}(\bm)$. Now, recall a typical situation of the mean-field theory with a first-order phase transition: There is a $J_\MF$ such that, for~$J$ near $J_\MF$, the mean-field free-energy function has two nearly degenerate minima separated by a barrier of height $\Delta(J)$, see Figure~\ref{fig1}. If the barrier~$\Delta(J)$ always exceeds the error term in~\eqref{5}, i.e., if $\Delta(J)>Jn\frac\kappa2 I_d$, some intermediate values of magnetization are forbidden and, as~$J$ increases through $J_\MF$, the physical magnetization undergoes a jump at some $J_{\text{\rm t}}$ near~$J_\MF$. See also Figure~\ref{fig2}.

\newcounter{obrazek}

\smallskip
The Main Theorem is a direct consequence of Theorem~\ref{thm1} and the following lemma:

\begin{Key Estimate}
Let $J\ge0$ and $\bb\in\E_\Omega$ and let $\bm_\star\in\mathscr{M}_\star(J,\bb)$. Let~$n$, $\kappa$ and $I_d$ be as in the Main Theorem. Then there is an infinite-volume Gibbs state $\nu_{J,\bb}$ for interaction \eqref{Ham} such that
\begin{equation}
\label{magav}
\bm_\star=\langle\bS_0\rangle_{J,\bb}
\end{equation} 
and
\begin{equation}
\label{cor-bd}
\bigl\langle(\bS_x,\bS_y)\bigr\rangle_{J,\bb}-|\bm_\star|^2\le n\kappa I_d,
\end{equation}
for any nearest-neighbor pair~$x,y\in\Z^d$. Here $\langle-\rangle_{J,\bb}$ denotes the expectation with respect to~$\nu_{J,\bb}$.
\end{Key Estimate}

The Key Estimate follows readily under certain conditions; for instance, when the parameter values~$J$ and~$\bb$ are such that there is a unique Gibbs state. Under these circumstances, the bound \eqref{cor-bd} is a special case of the infrared bound which can be derived using reflection positivity (see~\cite{FSS,DLS,FILS1,FILS2}) and paying close attention to the ``zero mode.'' Unfortunately, at the points of non-uniqueness, the bound in \eqref{cor-bd} is also needed. The restriction to extreme magnetizations is thus dictated by the need to approximate the magnetizations (and the states which exhibit them) by states where the standard ``RP, IRB'' technology can be employed.

The Key Estimate and Theorem~\ref{thm1} constitute a proof of the Main Theorem. Thus, a first-order phase transition (for $d\gg1$) can be  established in any  system of the form \eqref{Ham} by detailed analysis of the full mean-field theory. Although this sounds easy in principle, in practice there are cases where this can be quite a challenge. But, ultimately, the Main Theorem reduces the proof of a phase transitions to a problem in advanced calculus where (if desperate) one can employ computers to assist in the analysis.

\subsection{Direct argument for mean-field equation}
We have stated our main results in the context of the mean-field free energy. However, many practical calculations focus immediately on the mean-field equation for magnetization~\eqref{MFeq}. As it turns out, a direct study of the mean-field equation provides us with an alternative (albeit existential) approach to the results of this paper. The core of this approach is the variance bound for the magnetization stated as follows:

\begin{lemma}
\label{lemma1.3}
Let $d\ge 3$ and consider the spin system on $\Z^d$ with the Hamiltonian~\eqref{Ham}. Let~$n$ and $I_d$ be as in the Main Theorem. For $J\ge0$ and $\bb\in\E_\Omega$, let $\bm_\star\in\mathscr{M}_\star(J,\bb)$. Then there is an infinite-volume Gibbs state $\nu_{J,\bb}$ for the interaction \eqref{Ham} such that $\bm_\star=\langle\bS_0\rangle_{J,\bb}$ and
\begin{equation}
\label{magvar}
\Bigl\langle\Bigl|\frac1{2d}\sum_{x\colon |x|=1}\bS_x-\bm_\star\Bigr|^2\Bigr\rangle_{J,\bb}\le nJ^{-1}I_d,
\end{equation}
where $\langle-\rangle_{J,\bb}$ denotes the expectation with respect to $\nu_{J,\bb}$.
\end{lemma}

Here is how the bound \eqref{magvar} can be used to prove that mean-field equations are accurate in sufficiently large dimensions: Conditioning on the spin values at the neighbors of the origin and recalling the definition of $G(\bh)$, the expectation~$\langle\bS_0\rangle_{J,\bb}$ can be written as
\begin{equation}
\label{towardsMF}
\langle\bS_0\rangle_{J,\bb}=\biggl\langle\nabla G\biggl(\frac J{2d}\sum_{x\colon |x|=1}\bS_x+\bb\biggr)\biggr\rangle_{\!J,\bb}.
\end{equation}
Since the right-hand side of \eqref{magvar} tends to zero as $d\to\infty$, the (spatial) average of the spins neighboring the origin---namely $\frac 1{2d}\sum_{x\colon |x|=1}\bS_x$---is, with high probability, very close to $\bm_\star$. Using this in \eqref{towardsMF}, we thus find that $\bm_\star$ approximately satisfies the mean-field equation \eqref{MFeq}. Thus, to demonstrate phase coexistence (for $d\gg1$) it is sufficient to show that, along some curve in the parameter space, the solutions to the mean-field equations cannot be assembled into a continuous function. In many cases, this can be done dramatically by perturbative arguments.

While this alternative approach has practical appeal for certain systems, the principal drawback is  that  it provides no clue as to  the location of  the transition temperature. Indeed, as mentioned in the paragraph following the Main Theorem, secondary minima and other irrelevant  solutions to the mean-field equations typically develop well below $J=J_\MF$. Without the guidance of the free energy, there is no way of knowing which solutions are physically relevant. 

\section{Results for specific models}
\label{sec2}\medskip\noindent
In this section we adapt the previous general statements to three models: the~$q$-state Potts model, the $r$-component cubic model and the $O(N)$-nematic liquid crystal model. For appropriate ranges of the parameters~$q$, $r$ and~$N$ and dimension sufficiently large, we show that these models undergo a first-order phase transition as~$J$ varies. The relevant results appear as Theorems~\ref{Pottsres},~\ref{Cubicres} and~\ref{nematicres}.

\subsection{Potts model}
\label{sec2.1}\noindent
The Potts model, introduced in \cite{Potts}, is usually described as having a discrete spin space with~$q$ states, $\sigma_x\in\{1,2,\dots,q\}$, with the (formal) Hamiltonian
\begin{equation}
\beta H=-\eusmJ\sum_{\langle x,y\rangle}\delta_{\sigma_x,\sigma_y}.
\end{equation}
Here~$\delta_{\sigma_x\sigma_y}$ is the usual Kronecker delta and $\eusmJ=\frac J{2d}$. To bring the interaction into the form of \eqref{Ham}, we use the so called \emph{tetrahedral}
representation, see~\cite{Wu}. In particular, we let $\Omega=\{\hatv_1,\dots,\hatv_q\}$, where~$\hatv_\alpha$ denote the vertices of a $(q-1)$-dimensional hypertetrahedron, i.e.,
$\hatv_\alpha\in\R^{q-1}$ with
\begin{equation}
\hatv_\alpha\cdot\hatv_\beta=\begin{cases}1,\qquad&\text{if
}\alpha=\beta,\\
-\frac1{q-1},\qquad&\text{otherwise.}
\end{cases}
\end{equation}
The inner product is proportional to the usual dot product in $\R^{q-1}$. Explicitly, if~$\bS_x\in\Omega$ corresponds to~$\sigma_x\in\{1,\dots,q\}$, then we have  
\begin{equation}
(\bS_x,\bS_y)=\frac{q-1}q\,\bS_x\cdot\bS_y=\delta_{\sigma_x,\sigma_y}-\frac1q.
\end{equation}
(The reason for this rescaling the dot product is to maintain coherence with existing treatments of the mean-field version of this model.) The \emph{a priori} measure~$\mu$ gives a uniform weight to all~$q$~states in~$\Omega$.

Let us summarize some of the existing rigorous results about the~$q$-state Potts model. The $q=2$ model is the Ising model, which in mean-field theory as well as real life has a continuous transition. It is believed that the Potts model has a discontinuous transition for all $d\ge3$ and $q\ge3$ (see, e.g.,~\cite{Wu}).  In any $d\ge2$, it was first proved in~\cite{Kotecky-Shlosman} that for~$q$ sufficiently large, the energy density has a region of forbidden values over which it must jump discontinuously as~$J$ increases. On the basis of FKG monotonicity properties, see~\cite{ACCN}, this easily implies that the magnetization is also discontinuous. Such results have been refined and improved; for instance in \cite{KLMR,LMMRS}, Pirogov-Sinai type expansions have been used to show that there is a single point of discontinuity outside of which all quantities are analytic. However, for $d\ge3$, the values of~$q$ for which these techniques work are ``astronomical,'' and, moreover, deteriorate exponentially with increasing dimension.

Let~$m_\star(J)$ and~$e_\star(J)$ denote the the actual magnetization and energy density, respectively. These quantities can be defined using one-sided derivatives of the physical free energy:
\begin{equation}
\label{mstar-estar}
m_\star(J)=\frac\partial{\partial b}F(J,b\hatv_1)\Bigr|_{b=0^+}
\quad\text{and}\quad
e_\star(J)=\frac\partial{\partial J'}F(J',0)\Bigr|_{J'=J^+},
\end{equation}
or, equivalently, by optimizing the expectations $\langle (\hatv_1,\bS_0)\rangle$, resp., $\frac12\langle (\bS_0,\bS_x)\rangle$, where ``$0$'' is the origin and~$x$ is its nearest neighbor, over all Gibbs states that are invariant under the symmetries of~$\Z^d$. Recalling the Fortuin-Kasteleyn representation~\cite{Fortuin-Kasteleyn,ACCN,Gr,Georgii-Haggstrom-Maes}, let $P_\infty(J)$ be the probability that, in the associated random cluster model with parameters $p=1-e^{-J/(2d)}$ and~$q$, the origin lies in an infinite cluster. Then $m_\star(J)$ and $P_\infty(J)$ are related by the equation
\begin{equation}
m_\star(J)=\frac{q-1}qP_\infty(J).
\end{equation}
As a consequence, the magnetization~$m_\star(J)$ is a non-decreasing and right-continuous function of~$J$. The energy density~$e_\star(J)$ is non-decreasing in~$J$ simply by concavity of the free energy. The availability of the graphical representation allows us to make general statements about the phase-structure of these systems. In particular, in any $d\ge2$ and for all~$q$ under consideration, there is a $J_{\text{\rm c}}=J_{\text{\rm c}}(q,d)\in(0,\infty)$ such that $m_\star(J)>0$ for $J>J_{\text{\rm c}}$ while $m_\star(J)=0$ for $J<J_{\text{\rm c}}$, see~\cite{ACCN,Gr}. Whenever $m_\star(J_{\text{\rm c}})>0$ (which, by the aforementioned results~\cite{Kotecky-Shlosman,KLMR,LMMRS}, is known for $q\gg1$), there are at least~$q+1$ distinct extremal, translation-invariant Gibbs states at $J=J_{\text{\rm c}}$.

The mean-field free energy for the model without external field is best written in terms of components of $\bm$: If $(x_1,\dots,x_q)$ is a probability vector, we express $\bm$ as 
\begin{equation}
\label{bmPotts} \bm=x_1\hatv_1+\dots+x_q\hatv_q.
\end{equation}
The interpretation of this relation is immediate:~$x_k$ corresponds to the proportion of spins in the \hbox{$k$-th} spin-state. In terms of the variables in \eqref{bmPotts}, the mean-field free-energy function is (to within a constant) given by
\begin{equation}
\label{MFPotts} \varPhi_J(\bm)=\sum_{k=1}^q
\bigl(-{\textstyle\frac J2} x_k^2+x_k\log x_k\bigr).
\end{equation}
In \eqref{MFPotts} we have for once and all set the external field~$\bb$ to zero and suppressed it from the notation.

It is well-known (see~\cite{Wu,KS} and also Lemma~\ref{lemma3.1} of the present paper) that, for each $q\ge3$, there is a $J_\MF\in(2,q)$ such that~$\varPhi_J$ has a unique global minimizer $\bm=\bzero$ for $J<J_\MF$, while for $J>J_\MF$, there are~$q$ global minimizers which are obtained by permutations of single $(x_1,\dots,x_q)$ with~$x_1>x_2=\dots=x_q$. To keep the correspondence with $m_\star(J)$, we define the scalar mean-field magnetization~$m_\MF(J)$ as the maximal \emph{Euclidean} norm of all global minimizers of the mean-field free energy~$\varPhi_J(\bm)$. (In this parametrization, the asymmetric global maxima will be given by~$x_1=\frac1q+m_\MF(J)$ and~$x_2=\dots=x_q=\frac1q-\frac1{q-1}m_\MF(J)$.) Then $m_\MF(J)$ is the maximal positive solution to the equation
\begin{equation}
\label{mJPotts} 
\frac q{q-1}\,m=\frac{e^{J\frac q{q-1}m}-1}{e^{J\frac q{q-1}m}+q-1}.
\end{equation}
In particular, $J\mapsto m_\MF(J)$ is non-decreasing. We note that the explicit values of the coupling constant~$J_\MF$ and the magnetization $m_{\text{\rm c}}=m_\MF(J_\MF)$ at the mean-field transition are known:
\begin{equation}
\label{JMF}
J_\MF=2\,\frac{q-1}{q-2}\,\log(q-1)\quad\text{and}\quad
m_{\text{\rm c}}=\frac{q-2}q,
\end{equation}
see e.g.~\cite{Wu}. Thus, the mean-field transition is first-order for all $q>2$.

\smallskip
Our main result about the Potts model is then as follows:

\begin{theorem}[Potts model]
\label{Pottsres}
Consider the~$q$-state Potts model on $\Z^d$ and let $m_\star(J)$ be its scalar magnetization.
For each $q\ge3$, there exists a $J_{\text{\rm t}}=J_{\text{\rm t}}(q,d)$ and
two numbers $\epsilon_1=\epsilon_1(d,J)>0$ and $\epsilon_2=\epsilon_2(d)>0$ satisfying
$\epsilon_1(d,J)\to0$, uniformly on finite intervals of~$J$, and $\epsilon_2(d)\to0$ as $d\to\infty$, such that the following holds:
\begin{equation}
\label{Potts-mag-lb}
m_\star(J)\le \epsilon_1\quad\text{for}\quad J<J_{\text{\rm t}}
\end{equation}
and
\begin{equation}
\label{Potts-mag-ub}
|m_\star(J)-m_\MF(J)|\le \epsilon_1\quad\text{for}\quad J>J_{\text{\rm t}}.
\end{equation}
Moreover,
\begin{equation}
\label{Potts-J-bd}
|J_{\text{\rm t}}-J_\MF|\le \epsilon_2.
\end{equation}
In particular, both the magnetization $m_\star(J)$ and the energy density $e_\star(J)$ undergo a jump at $J=J_{\text{\rm t}}$ whenever $d$ is sufficiently large.
\end{theorem}

The jump in the energy density at~$J_{\text{\rm t}}$ immediately implies the existence of at least $q+1$ distinct extremal Gibbs measures at $J=J_{\text{\rm t}}$. However, the nature of our proofs does not permit us to conclude that $m_\star(J)=0$ for $J<J_{\text{\rm t}}$ nor can we rule out that~$m_\star(J)$ undergoes further jumps for $J>J_{\text{\rm t}}$. (Nonetheless, the jumps for $J>J_{\text{\rm t}}$ would have to be smaller than $2\epsilon_1(d)$.) Unfortunately, we can say nothing about the continuous-$q$ variant of the Potts model---the random cluster model---for non-integer~$q$. In this work, the proofs lean too heavily on the spin representation. Furthermore, for non-integer~$q$, the use of our principal tool, reflection positivity, is forbidden;~see~\cite{Biskup}.

We also concede that, despite physical intuition to the contrary, our best bounds on~$\epsilon_2(d)$ and~$\epsilon_1(d,J)$ deteriorate with increasing~$q$. This is an artifact of the occurrence of the single-spin space dimension on the right-hand side of \eqref{5}. (This sort of thing seems to plague all existing estimates based on reflection positivity.) In particular, we cannot yet produce a sufficiently large dimension $d$ for which the phase transition in all ($q\ge3$)-state Potts models would be provably~first~order. 

\subsection{Cubic model{}}
\label{sec2.2}\noindent
Our second example of interest is the $r$-component cubic model. 
Here the spins~$\bS_x$ are the unit vectors in the coordinate directions of~$\R^r$, i.e., if $\hate_k$ are the standard unit vectors
in~$\R^r$, then
\begin{equation}
\Omega=\{\pm\hate_k\colon k=1,\dots, r\}.
\end{equation}
The Hamiltonian is given by \eqref{Ham}, with the inner product given by the usual dot product in~$\R^r$ and the \emph{a priori} measure given by the uniform measure on~$\Omega$. As in the last subsection, we set $\bb=\bzero$ and suppress any $\bb$-dependence from the notation. We note that the $r=1$ case is the Ising model while the case $r=2$ is equivalent to two uncoupled Ising models.

The cubic model was introduced (and studied) in \cite{FUffer1,FUffer2} 
as a model of the magnetism in rare-earth compounds with a cubic crystal symmetry. There it was noted that the associated mean-field theory has a discontinuous transition for $r\ge4$, while the transition is continuous for $r=1$,~$2$ and~$3$. The mean field theory is best expressed in terms of the collection of parameters $\bar y=(y_1,\dots,y_r)$ and $\bar\mu=(\mu_1,\dots,\mu_r)$, where $y_k$ stands for the fraction of spins that take the values $\pm\hate_k$ and $\mu_ky_k$ is the magnetization in the direction~$\hate_k$. In this language, the magnetization vector can be written as 
\begin{equation}
\label{bmrel}
\bm=y_1\mu_1\hate_1+\dots+y_r\mu_r\hate_r.
\end{equation}
To describe the mean-field free-energy function, we define
\begin{equation}
\label{MFcubic} 
K^{(r)}_J(\bar y,\bar\mu)=\sum_{k=1}^r\bigl(y_k\log y_k+
y_k\,\varTheta_{2Jy_k}\!(\mu_k)\bigr),
\end{equation}
where $\varTheta_J(\mu)$ denotes the standard Ising mean-field free energy with bias~$\mu$; i.e., the quantity in~\eqref{MFPotts} with $q=2$,~$x_1=\frac12(1+\mu)$ and~$x_2=\frac12(1-\mu)$. Then~$\varPhi_J(\bm)$ is found by minimizing $K^{(r)}_J(\bar y,\bar\mu)$ over all allowed pairs $(\bar y,\bar\mu)$ such that \eqref{bmrel} holds.

As in the case of the Potts model, the global minimizer of~$\varPhi_J(\bm)$ will be a permutation of a highly-symmetric state. However, this time the result is not so well known, so we state it as a separate proposition:

\begin{proposition}
\label{prop-cubic}
Consider the $r$-component cubic model. For each $J\ge 0$, the only local minima of~$\varPhi_J$ are $\bm=\bzero$ or~$\bm=\pm m_\MF\,\hate_k$, $k=1,\dots,r$, where $m_\MF=m_\MF(J)$ is the maximal positive solution to the equation
\begin{equation}
\label{cubic-eq}
m=\frac{\sinh Jm}{r-1+\cosh Jm}.
\end{equation}
Furthermore, there is a $J_\MF\in(0,\infty)$ such that the only global minimizers of~$\varPhi_J(\bm)$ are $\bm=\bzero$ for $J<J_\MF$ and~$\bm=\pm m_\MF(J)\hate_k$, $k=1,\dots,r$, (with $m_\MF(J)>0$) for $J>J_\MF$.
\end{proposition}

For a system on $\Z^d$, the scalar magnetization is most conveniently defined as the norm of $\langle\bS_0\rangle_J$, optimized over all translation-invariant Gibbs states for the coupling constant~$J$. The energy density~$e_\star(J)$ is defined using the same formula as for the Potts model, see~\eqref{mstar-estar}.

\smallskip
Our main result about the cubic model is then as follows:

\begin{theorem}[Cubic model]
\label{Cubicres}
Consider the $r$-state cubic model on $\Z^d$ and let $m_\star(J)$ be its scalar magnetization. Then for every $r\ge4$, there exists a $J_{\text{\rm t}}=J_{\text{\rm t}}(q,d)$ and
two numbers $\epsilon_1=\epsilon_1(d,J)>0$ and $\epsilon_2=\epsilon_2(d)>0$ satisfying $\epsilon_1(d,J)\to0$, uniformly on finite intervals of~$J$, and $\epsilon_2(d)\to0$ as $d\to\infty$, such that the following holds:
\begin{equation}
\label{cubic-down}
m_\star(J)\le \epsilon_1\quad\text{for}\quad J<J_{\text{\rm t}}
\end{equation}
and
\begin{equation}
\label{cubic-up}
|m_\star(J)-m_\MF(J)|\le \epsilon_1\quad\text{for}\quad J>J_{\text{\rm t}}.
\end{equation}
Moreover,
\begin{equation}
\label{cubic-J-bd}
|J_{\text{\rm t}}-J_\MF|\le \epsilon_2.
\end{equation}
In particular, both the magnetization $m_\star(J)$ and the energy density $e_\star(J)$ undergo a jump at $J=J_{\text{\rm t}}$ whenever $d$ is sufficiently large.
\end{theorem}

As in the case of the Potts model, our technique does not allow us to conclude that $J_{\text{\rm t}}$ is the only value of~$J$ where the magnetization undergoes a jump. In this case, we do not even know that the magnetization is a monotone function of~$J$; the conclusions \twoeqref{cubic-down}{cubic-up} can be made because we know that the energy density is close to~$\frac12m_\star(J)^2$ and is (as always) a non-decreasing function of~$J$. Finally, we also cannot prove that, in the state with large magnetization in the direction~$\hate_1$, there will be no additional symmetry breaking in the other directions. Further analysis, based perhaps on graphical representations, is needed.

\subsection{Nematic liquid-crystal model}
\label{sec2.3}\noindent
The nematic models are designed to study the behavior of liquid crystals, see the monograph~\cite{deGennes} for more background on the subject. In the simplest cases, 
a liquid crystal may be regarded as a suspension of rod-like molecules which, for all intents and purposes, are symmetric around their midpoint. For the models of direct physical relevance, each rod (or a small collection of rods) is described by an three-dimensional spin and one considers only interactions that are (globally) $O(3)$-invariant and invariant under the (local) reversal of any spin. The simplest latticized version of such a system is described by the Hamiltonian
\begin{equation}
\label{HamNM}
\beta H(\bs)=-\frac J{2d}\sum_{\langle x,y\rangle}(\bs_x\cdot\bs_y)^2,
\end{equation}
with~$\bs_x$ a unit vector in~$\R^3$ and~$x\in\Z^d$ with~$d=2$ or~$d=3$. We will study the above Hamiltonian, but we will consider general dimensions~$d$ (provided $d\ge3$) and spins that are unit vectors in any~$\R^N$ (provided $N\ge3$).

The Hamiltonian \eqref{HamNM} can be rewritten into the form \eqref{Ham} as follows~\cite{deGennes}: Let $\E_\Omega$ be the space of all traceless $N\times N$ matrices with real coefficients and let $\Omega$ be the set of those matrices $\bQ=(Q_{\alpha,\beta})\in\E_\Omega$ for which there is a unit vector in $\bv=(v_\alpha)\in\R^N$ such that
\begin{equation}
\label{1.30}
Q_{\alpha\beta}=v_\alpha v_\beta-\frac1N\delta_{\alpha\beta},\qquad \alpha,\beta=1,\dots,N.
\end{equation}
Writing $\bQ_x$ for the matrix arising from the spin $\bs_x$ via \eqref{1.30}, the interaction term becomes
\begin{equation}
\label{1.31}
(\bs_x\cdot\bs_y)^2=\text{Tr}(\bQ_x\bQ_y)+\frac1N.
\end{equation}
Now $\E_\Omega$ is a finite-dimensional vector space and $(\bQ,\bQ')=\text{Tr}(\bQ\bQ')$ is an inner product on $\E_\Omega$, so \eqref{HamNM} indeed takes the desired form \eqref{Ham}, up to a constant that has no relevance for physics.

The \emph{a priori} measure on $\Omega$ is a pull-back of the uniform distribution on the unit sphere in $\R^N$. More precisely, if $\bv$ is uniformly distributed on the unit sphere in $\R^N$, then $\bQ\in\Omega$ is a random variable arising from $\bv$ via \eqref{1.30}. As a consequence, the \emph{a priori} distribution is invariant under the action of the Lee group $O(N,\R)$ given by
\begin{equation}
\label{gaction}
\bQ_x\mapsto \bg^{-1}\bQ_x\bg,\qquad \bg\in O(N,\R).
\end{equation}
The parameter signaling the phase transition, the so called \emph{order parameter}, is ``tensor'' valued. In particular, it corresponds to the expectation of~$\bQ_0$. The order parameter can always be diagonalized. The diagonal form is not unique; however, we can find an orthogonal transformation that puts the eigenvalues in a decreasing order. Thus the order parameter is effectively an~$N$-vector $\blambda=(\lambda_1,\dots,\lambda_N)$ such that $\lambda_1\ge\lambda_2\ge\dots\ge\lambda_N$. We note that, since each $\bQ_x$ is traceless, we~have~$\sum_k\lambda_k=0$.

The previous discussion suggests the following definition of the \emph{scalar} order parameter: For $J\ge0$, we let $\lambda_\star(J)$ be the value of the largest non-negative eigenvalue of the matrix $\langle \bQ_0\rangle_J$, optimized over all translation-invariant Gibbs states for the coupling constant~$J$. As far as rigorous results about the quantity $\lambda_\star(J)$ are concerned, we know from~\cite{Angelescu-Zagrebnov} that (in $d\ge3$) $\lambda_\star(J)>0$ once~$J$ is sufficiently large. On the other hand, standard high-temperature techniques (see 
e.g.~\cite{Dobrushin,vdBerg-Maes,Alexander-Chayes}) show that if~$J$ is sufficiently small then there is a unique Gibbs state. In particular, since this state is then invariant under the action \eqref{gaction} of the full $O(N,\R)$ group, this necessitates that $\lambda_\star(J)\equiv0$ for~$J$ small enough. The goal of this section is to show that~$\lambda_\star(J)$ actually undergoes a \emph{jump} as~$J$ varies.

The mean-field theory of the nematic model is formidable. Indeed, for any particular~$N$  it does not seem possible to obtain a workable expression for~$\varPhi_J(\blambda)$, even if we allow that the components of~$\blambda$ have only two distinct values (which is usually assumed without apology in the physics literature). Notwithstanding, this simple form of the vector minimizer and at least some of the anticipated properties can be established:

\begin{proposition}
\label{prop2.4a}
Consider the $O(N)$-nematic model for $N\ge3$. Then every local minimum of~$\varPhi_J(\blambda)$ is an orthogonal transformation of the matrix $\blambda = \diag(\lambda,-\frac\lambda{N-1},\dots,-\frac\lambda{N-1})$,
where $\lambda$ is a non-negative solution to the equation
\begin{equation}
\label{lMFeq}
\lambda=\frac{\displaystyle\int_0^1\textd x\,(1-x^2)^{\frac{N-3}2}
e^{\frac{JN\lambda}{N-1}x^2}\bigl(x^2-\tfrac1N\bigr)}
{\displaystyle\int_0^1\textd x\,(1-x^2)^{\frac{N-3}2}
e^{\frac{JN\lambda}{N-1}x^2}}.
\end{equation}
In particular, there is an increasing and right-continuous function $J\mapsto\lambda_\MF(J)$ such that the unique minimizer of~$\varPhi_J(\blambda)$ is $\blambda=\bzero$ for $J<J_\MF$, while for any $J>J_\MF$, the function~$\varPhi_J(\blambda)$ is minimized by the orthogonal transformations of
\begin{equation}
\label{diag-nematics}
\blambda=\diag\Bigl(\lambda_\MF(J), -\frac{\lambda_\MF(J)}{N-1},\dots,-\frac{\lambda_\MF(J)}{N-1}\Bigr).
\end{equation}
At the continuity points of $\lambda_\MF\colon(J_\MF,\infty)\to[0,1]$, these are the only global minimizers of~$\varPhi_J$.
\end{proposition}

Based on the pictorial solution of the problem by physicists, see e.g.~\cite{deGennes}, we would expect that $J\mapsto\lambda_\MF(J)$ is continuous on its domain and, in fact, corresponds to the maximal positive solution to \eqref{lMFeq}. (This boils down to showing certain convexity-concavity property of the function on the right-hand side of \eqref{lMFeq}.) While we could not establish this fact for all $N\ge3$, we were  successful at least for~$N$ sufficiently large. The results of the large-$N$ analysis are summarized as~follows:

\begin{proposition}
\label{prop2.4b}
Consider the $O(N)$-nematic model for $N\ge3$ and let $\lambda_\MF^{(N)}(J)$ be the maximal positive solution to \eqref{lMFeq}. Then there exists an $N_0\ge3$ and, for each $N\ge N_0$, a number $J_\MF=J_\MF(N)\in(0,\infty)$ such that for each $N\ge N_0$, the unique minimizer of~$\varPhi_J(\blambda)$ is $\blambda=\bzero$ for $J<J_\MF$, while for any $J>J_\MF$, the function~$\varPhi_J(\blambda)$ is minimized only by the orthogonal transformations~of~\eqref{diag-nematics}, with $\lambda_\MF(J)>0$.

The function $J\mapsto\lambda_\MF^{(N)}(J)$ is continuous and strictly increasing on its domain and has the following large-$N$ asymptotic: For all $J\ge2$,
\begin{equation}
\label{limlam}
\lim_{N\to\infty}\lambda_\MF^{(N)}(JN)=\frac12\bigl(1+\sqrt{1-4J^{-2}}\bigr).
\end{equation}
Moreover, there exists a $J_\MF^{(\infty)}$ (with $J_\MF^{(\infty)}\approx2.455$) such that
\begin{equation}
\label{limJ}
\lim_{N\to\infty} \frac{J_\MF(N)}N=J_\MF^{(\infty)}.
\end{equation}
\end{proposition}

Now we are ready to state our main theorem concerning $O(N)$-nematics.
As can be gleaned from a \emph{careful} reading, our conclusions are not quite as strong as in the previous cases
(due the intractability of the associated mean-field theory).
Nevertheless, a \emph{bona fide} first-order transition is established for these systems.

\begin{theorem}[Nematic model]
\label{nematicres}
Consider the $O(N)$-nematic model with the Hamiltonian \eqref{HamNM} and $J\ge0$. For each $N\ge3$, there exists a non-negative function $J\mapsto\lambda_\MF^\star(J)$, a constant $J_{\text{\rm t}}=J_{\text{\rm t}}(N,d)$ and two numbers $\epsilon_1=\epsilon_1(d,J)>0$ and $\epsilon_2=\epsilon_2(d)>0$ satisfying $\epsilon_1(d,J)\to0$, uniformly on finite intervals of~$J$, and $\epsilon_2(d)\to0$ as $d\to\infty$, such that the following holds: 

For all $J\ge0$, the matrix $\blambda=\diag(\lambda_\MF^\star(J),-\frac{\lambda_\MF^\star(J)}{N-1}, \dots,-\frac{\lambda_\MF^\star(J)}{N-1})$ is a local minimum of~$\varPhi_J$. Moreover, we have the bounds
\begin{equation}
\label{nematic-down}
\lambda_\star(J)\le\epsilon_1\quad\text{for}\quad J<J_{\text{\rm t}}
\end{equation}
and
\begin{equation}
\label{nematic-up}
|\lambda_\star(J)-\lambda_\MF^\star(J)|\le\epsilon_1
\quad\text{for}\quad J>J_{\text{\rm t}}.
\end{equation}
Furthermore,
\begin{equation}
|J_{\text{\rm t}}-J_\MF|\le \epsilon_2.
\end{equation}
In particular, $\lambda_\star(J)\ge\varkappa>0$ for all $J>J_{\text{\rm t}}$ and all $N\ge3$ and both the order parameter and the energy density $e_\star(J)$ undergo a jump at $J=J_{\text{\rm t}}$, provided the dimension is sufficiently large. 
\end{theorem}

The upshot of the previous theorem is that the high-temperature region with $\blambda=\bzero$ and the low-temperature region with $\blambda\ne\bzero$ (whose existence was proved in \cite{Angelescu-Zagrebnov}) are separated by a first-order transition. However, as with the other models, our techniques are not sufficient to prove that~$\blambda$ is exactly zero for all $J<J_{\text{\rm t}}$, nor, for $J>J_{\text{\rm t}}$, that all states are devoid of some other additional breakdown of symmetry. Notwithstanding, general theorems about Gibbs measures guarantee that, a jump of $J\mapsto\lambda_\star(J)$ at $J=J_{\text{\rm t}}$ implies the coexistence of a ``high-temperature'' state with various symmetry-broken ``low-temperature'' states.

\section{Proofs of mean-field bounds}
\label{sec3}\noindent
\vspace{-6mm} 
\subsection{Convexity estimates}
\label{sec3.1}\noindent
In order to prove Theorem~\ref{thm1}, we need to recall a few standard notions from convexity  theory and prove a simple lemma. Let $\mathscr{A}\subset\R^n$ be a convex set. Then we  define  the \emph{affine hull} of~$\mathscr{A}$ by the formula
\begin{equation}
\aff\mathscr{A}=\bigl\{\lambda x+(1-\lambda)y\colon x,y\in\mathscr{A},\,\lambda\in\R\bigr\}.
\end{equation}
(Alternatively, $\aff\mathscr{A}$ is an smallest affine subset of~$\R^n$ containing~$\mathscr{A}$.) This concept allows us to define the \emph{relative interior}, $\ri\mathscr{A}$, of~$\mathscr{A}$ as the set of all~$x\in\mathscr{A}$ for which there exists an $\epsilon>0$ such that
\begin{equation}
\label{ricond}
y\in\aff\mathscr{A}\quad\&\quad|y-x|\le\epsilon\quad\Rightarrow\quad y\in\mathscr{A}.
\end{equation}
It is noted that this definition of relative interior differs from the standard  topological definition. For us it is important that the standard (topological) closure of $\ri\mathscr{A}$ is simply the standard  closure of~$\mathscr{A}$. We refer to \cite{Rocky} for more details.

\begin{lemma}
\label{lemma3.0}
For each $\bm\in\ri\{\bm'\in\E_\Omega\colon S(\bm')>-\infty\}$, there exists a vector $\bh\in\E_\Omega$ such that \hbox{$\nabla G(\bh)=\bm$}.
\end{lemma}

Results of this sort are quite well known; e.g., with some effort this can be gleaned from Lemma~2.2.12 in~\cite{Dembo-Zeitouni} combined with the fact that the so called exposed points of $S(\bm)$ can be realized as $\nabla G(\bh)$ for some~$\bh$. For completeness, we provide a full derivation which exploits the particulars of the setup at hand.
\smallskip

\begin{proofsect}{Proof}
Let~$\mathscr{C}$ abbreviate $\{\bm'\in\E_\Omega\colon S(\bm')>-\infty\}$ and let $\bm\in\ri\mathscr{C}$. Let us define the set $\V=\{\bm'-\bm\colon \bm'\in\aff\mathscr{C}\}$. It is easy to see that~$\V$ is in fact the affine hull of the  shifted set  $\mathscr{C}-\bm$  and, since $\bzero\in\V$, it is a closed linear subspace of $\E_\Omega$. First we 
claim that the infimum in \eqref{S} can be restricted to $\bh\in\V$. Indeed, if $\bh,\ba\in\E_\Omega$, then the convexity of $\bh\mapsto G(\bh)$ gives
\begin{equation}
\label{cG}
G(\bh+\ba)-(\bh+\ba,\bm)\ge G(\bh)-(\bh,\bm) + \bigl(\ba,\nabla G(\bh)-\bm\bigr)
\end{equation}
for any $\bm$. This implies that $\nabla G(\bh)$ has a finite entropy, i.e., $\nabla G(\bh)\in\mathscr{C}$ for any $\bh\in\E_\Omega$. Now let~$\bm$ be as above and $\ba\in\V^\bot$. Then an inspection of the definition of~$\V$ shows that the last term in \eqref{cG} identically vanishes. Consequently, for the infimum \eqref{S}, we will always be better off  with $\bh\in\V$.

Let $\bh_k\in\V$ be a minimizing sequence for $S(\bm)$; i.e., $G(\bh_k)-(\bh_k,\bm)\to S(\bm)$ as $k\to\infty$. We claim that $\bh_k$ contains a subsequence tending to a finite limit. Indeed, if on the contrary $h_k=|\bh_k|\to\infty$ we let $\btau_k$ be defined by $\bh_k=h_k\btau_k$ and suppose that $\btau_k\to\btau$ (at least along a subsequence), where $|\btau|=1$. Now since $\bm\in\ri\mathscr{C}$ and $\btau\in\V$, we have $\bm+\epsilon\btau\in\aff\mathscr{C}$ for all $\epsilon$ and, by \eqref{ricond}, $\bm+\epsilon\btau\in\mathscr{C}$ for some $\epsilon>0$ sufficiently small. But we also have
\begin{equation}
G(\bh_k)-(\bh_k,\bm+\epsilon\btau)=G(\bh_k)-(\bh_k,\bm)-\epsilon h_k(\btau_k,\btau),
\end{equation}
which tends to the negative infinity because $(\btau_k,\btau)\to1$ and $h_k\to\infty$. But then $S(\bm+\epsilon\btau)=-\infty$, which contradicts that $\bm+\epsilon\btau\in\mathscr{C}$. Thus $\bh_k$ contains a converging subsequence, $\bh_{k_j}\to\bh$. Using that $\bh$ is an actual minimizer of $G(\bh)-(\bh,\bm)$, it follows that $\nabla G(\bh)=\bm$.
\end{proofsect}

Now we are ready to prove our principal convexity bound:

\begin{proofsect}{Proof of Theorem~\ref{thm1}}
Recall that $F_\MF(J,\bb)$ denotes the infimum of~$\varPhi_{J,\bb}(\bm)$ over all $\bm\in\conv(\Omega)$.
As a first step, we will prove that there is a constant $C<\infty$ such that for any finite $\Lambda\subset\Z^d$ and any boundary condition $\bS'_{\partial\Lambda}$, the partition function obeys the bound
\begin{equation}
\label{Zlow}
Z_\Lambda(\bS'_{\partial\Lambda})\ge e^{-|\Lambda|F_\MF(J,\bb)-C|\partial\Lambda|},
\end{equation}
where $|\Lambda|$ denotes the number of sites in $\Lambda$ and $|\partial\Lambda|$ denotes the number of bonds of $\Z^d$ with one end in $\Lambda$ and the other in $\Z^d\setminus\Lambda$. (This is an explicit form of the well known fact that the free energy is always lower than the associated mean-field free energy, see~\cite{Ellis,Simon}.)

To prove \eqref{Zlow}, let $\bM_\Lambda$ denote the total magnetization in $\Lambda$,
\begin{equation}
\label{MLambda}
\bM_\Lambda=\sum_{x\in\Lambda}\bS_x,
\end{equation}
and let $\langle-\rangle_{0,\bh}^{(\Lambda)}$ be the \emph{a priori} state in~$\Lambda$ tilted with a uniform magnetic field $\bh$, i.e., for any measurable function~$f$ of the configurations in~$\Lambda$,
\begin{equation}
\label{3.3}
\langle f\rangle_{0,\bh}^{(\Lambda)}=e^{-|\Lambda|G(\bh)}\langle f e^{(\bh,\bM_\Lambda)}\rangle_0.
\end{equation}
Fix an $\bh\in\E_\Omega$ and let $\bm_\bh=\nabla G(\bh)$. By inspection, $\nabla G(\bh)=\langle\bS_x\rangle_{0,\bh}^{(\Lambda)}$ for all~$x\in\Lambda$. Then 
\begin{equation}
Z_\Lambda(\bS'_{\partial\Lambda})=e^{|\Lambda|G(\bh)}\bigl
\langle e^{-(\bh,\bM_\Lambda)-\beta H_\Lambda(\bS_\Lambda|\bS'_{\partial\Lambda})}\bigr\rangle_{0,\bh}^{(\Lambda)},
\end{equation}
which using Jensen's inequality gives
\begin{equation}
\label{inter}
Z_\Lambda(\bS'_{\partial\Lambda})\ge\exp\Bigl\{|\Lambda|\bigl(G(\bh)-(\bh,\bm_\bh)\bigr)-\bigl\langle\beta H(\bS_\Lambda|\bS'_{\partial\Lambda})\bigr\rangle_{0,\bh}\Bigr\}.
\end{equation}
To estimate the expectation of $\beta H(\bS_\Lambda|\bS'_{\partial\Lambda})$, we first discard (through a bound) the boundary terms and then evaluate the contribution of the interior bonds.
Since the number of interior bonds in $\Lambda$ is more than $d|\Lambda|-|\partial\Lambda|$, this gets us
\begin{equation}
-\bigl\langle\beta H(\bS_\Lambda|\bS'_{\partial\Lambda})\bigr\rangle_{0,\bh}\ge
-\frac J2|\bm_\bh|^2-C|\partial\Lambda|.
\end{equation}
Now $G(\bh)-(\bh,\bm_\bh)\ge S(\bm_\bh)$, so we have $Z_\Lambda(\bS'_{\partial\Lambda})\ge e^{-|\Lambda|\varPhi_{J,\bb}(\bm_\bh)-C|\partial\Lambda|}$.  But Lemma~\ref{lemma3.0} guarantees that each $\bm$ with $S(\bm)>-\infty$ can be approximated by a sequence of $\bm_\bh$ with $\bh\in\E_\Omega$, so the bound \eqref{Zlow} follows by optimizing over $\bh\in\E_\Omega$.

Next, let $\nu_{J,\bb}$ be an infinite volume Gibbs state and let $\langle-\rangle_{J,\bb}$ denote expectation with respect to $\nu_{J,\bb}$. Then we claim that
\begin{equation}
\label{2.6}
e^{|\Lambda|G(\bh)}=\bigl\langle e^{(\bh,\bM_\Lambda)+\beta H_\Lambda(\bS_\Lambda|\bS_{\partial\Lambda})}Z_\Lambda(\bS_{\partial\Lambda})\bigr\rangle_{J,\bb}.
\end{equation}
(Here $\bS_\Lambda$, resp. $\bS_{\partial\Lambda}$ denote the part of the \emph{same} configuration $\bS$ inside, resp., outside~$\Lambda$. Note that the relation looks trivial for $\bh=0$.) Indeed, the conditional distribution in $\nu_{J,\bb}$ given that the configuration outside $\Lambda$ equals $\bS'$ is $\nu^{(\bS')}_\Lambda$, as defined in \eqref{nuSS}. But then \eqref{nuSS} tells us that
\begin{equation}
\int e^{(\bh,\bM_\Lambda)+\beta H_\Lambda(\bS_\Lambda|\bS')}Z_\Lambda(\bS')
\,\nu^{(\bS')}_\Lambda(\textd\bS_\Lambda)
=\int e^{(\bh,\bM_\Lambda)}\mu(\textd\bS_\Lambda)=e^{|\Lambda|G(\bh)}.
\end{equation}
The expectation over the boundary condition $\bS'$ then becomes irrelevant and \eqref{2.6} is proved.

Now suppose that $\nu_{J,\bb}$ is the $\Z^d$-translation and rotation invariant Gibbs measure in question and recall that $\bm_\star=\langle \bS_0\rangle_{J,\bb}$, where $\langle-\rangle_{J,\bb}$ denotes the expectation with respect to~$\nu_{J,\bb}$. To prove our desired estimate, we use \eqref{Zlow} on the right-hand side of \eqref{2.6} and apply Jensen's inequality to get
\begin{equation}
\label{2.8}
e^{|\Lambda|G(\bh)}\ge\exp\Bigl\{\bigl\langle(\bh,\bM_\Lambda)+\beta H_\Lambda\bigr\rangle_{J,\bb}\Bigr\}e^{-|\Lambda|F_\MF(J,\bb)-C|\partial\Lambda|}.
\end{equation}
Using the invariance of the state $\nu_{J,\bb}$ with respect to the translations and rotations of $\Z^d$, we have 
\begin{equation}
\bigl\langle(\bh,\bM_\Lambda)\bigr
\rangle_{J,\bb}=|\Lambda|(\bh,\bm_\star)
\end{equation}
while 
\begin{equation}
\langle\beta H_\Lambda\rangle_{J,\bb}\ge-|\Lambda|\frac J2\bigl\langle(\bS_0,\bS_x)\bigr\rangle_{J,\bb}
-|\Lambda|(\bb,\bm_\star)-C'|\partial\Lambda|,
\end{equation}
where~$C'$ is a constant that bounds the worst-case boundary term
and where~$x$ stands for any neighbor of the origin.
By plugging these bounds back into \eqref{2.8} and passing to the thermodynamic limit, we conclude that
\begin{equation}
-G(\bh)+(\bh-\bb,\bm_\star)-\frac J2\bigl\langle(\bS_0,\bS_x)
\bigr\rangle_{J,\bb}\le F_\MF(J,\bb).
\end{equation}
Now optimizing the left-hand side over $\bh\in\E_\Omega$ allows us to replace $-G(\bh)+(\bh,\bm_\star)$ by $-S(\bm_\star)$. Then the bound \eqref{genbd} follows by adding and subtracting the term $\frac J2|\bm_\star|^2$ on the left-hand side.
\end{proofsect}

\subsection{Infrared bound}
Our proof of the Key Estimate (and hence the Main Theorem) requires the use of the \emph{infrared bounds}, which in turn are derived from reflection positivity. The connection between infrared bounds and reflection positivity dates back (at least) to \cite{FSS,DLS,FILS1,FILS2}. However, the present formulation (essentially already contained in~\cite{FSS,BKLS,KS}) emphasizes more explicitly the role of the ``$k=0$'' Fourier mode of the two-point correlation function by subtracting the square of the background average.

Reflection positivity is greatly facilitated by first considering finite systems with periodic boundary conditions. If it happens that there is a \emph{unique} Gibbs state for parameter values~$J$ and~$\bb$ then the proof of the Key Estimate is straightforward---there is no difficulty with putting the system on a torus and taking the limit. In particular, the Key Estimate amounts (more or less) to Corollary~2.5 in~\cite{FSS}. But when there are several infinite-volume Gibbs states, we can anticipate trouble with the naive limits of the finite-volume torus states. Fortunately, Gibbsian uniqueness is not essential to our arguments. Below we list two properties of Gibbs states which allow a straightforward proof of the desired infrared bound. Then we show that in general we can obtain the infrared bound for states of interest by an approximation argument.

\begin{property}
\label{def1}
An infinite-volume Gibbs measure $\nu_{J,\bb}$\/\/ (not necessarily extremal) for the interaction \eqref{Ham} is called a {\rm torus state} if it can be obtained by a (possibly subsequential) weak limit as $L\to\infty$ of the Gibbs states in volume $[-L,L]^d\cap\Z^d$, for the interaction \eqref{Ham} with periodic boundary conditions.
\end{property}

Given~$J$ and~$\bb$, we let $\mathscr{M}(J,\bb)$ denote the subset of $\conv(\Omega)$ containing all magnetizations achieved by infinite-volume translation-invariant Gibbs states for the interaction \eqref{Ham}. Next, recall the notation~$\bM_\Lambda$ from \eqref{MLambda} for the average magnetization in $\Lambda\subset\Z^d$. 

\begin{property}
\label{def2}
An infinite-volume Gibbs measure $\nu_{J,\bb}$\/\/ (not necessarily extremal) for the interaction \eqref{Ham} is said to have {\rm block-average magnetization} $\bm$ if
\begin{equation}
\lim_{\Lambda\nearrow\Z^d}\frac{\bM_\Lambda}{|\Lambda|}=\bm,\qquad
\nu_{J,\bb}\text{-almost surely}.
\end{equation}
Here the convergence $\Lambda\nearrow\Z^d$ is along the net of all the finite boxes~$\Lambda\subset\Z^d$ with partial order induced by set inclusion. (See~\cite{Georgii} for more details.)
\end{property}

Our first goal is to show that every torus state with a deterministic block-average magnetization satisfies the infrared bound. Suppose $d\ge3$ and let $D^{-1}$ denote the Fourier transform of the inverse lattice Laplacian with Dirichlet boundary condition. In lattice coordinates, $D^{-1}$ has the representation
\begin{equation}
D^{-1}(x,y)=\int_{[-\pi,\pi]^d} \frac{\textd^d
k}{(2\pi)^d}\frac1{\widehat D(k)}\,e^{ik(x-y)}, \qquad x,y\in\Z^d,
\end{equation}
where $\widehat D(k)=1-\frac1d\sum_{j=1}^d\cos(k_j)$. Note that the integral converges by our assumption that $d\ge3$.

\begin{lemma}
\label{lemma2.2}
Let $d\ge3$ and suppose that $\nu_{J,\bb}$ is a Gibbs state for interaction \eqref{Ham} satisfying Properties~1 and~2. Let $\langle-\rangle_{J,\bb}$ denote the expectation with respect to $\nu_{J,\bb}$ and let~$\bm$ denote the value of magnetization in~$\nu_{J,\bb}$. Then for all $(v_x)_{x\in\Z^d}$ such that $v_x\in\R$ and $\sum_{x\in\Z^d}|v_x|<\infty$,
\begin{equation}
\label{IRBrel}
\sum_{x,y\in\Z^d}v_xv_y\,\bigl\langle(\bS_x-\bm,
\bS_y-\bm)\bigr\rangle_{J,\bb}\le nJ^{-1}
\sum_{x,y\in\Z^d}v_xv_y \,D^{-1}(x,y).
\end{equation}
Here~$n$ denotes the dimension of $\E_\Omega$.
\end{lemma}

\begin{proofsect}{Proof}
Let $\Lambda_L=[-L,L]^d\cap\Z^d$ and let $\nu_{J,\bb}^{(L)}$ be the finite-volume Gibbs state in $\Lambda_L$ for the interaction \eqref{Ham} with periodic boundary conditions. Let $\Lambda_L^\star=\{(\frac{2\pi}{2L+1}n_1,\dots,\frac{2\pi}{2L+1} n_d)\colon -L\le n_i\le L\}$ denote the reciprocal lattice. Let $(\bw_x)_{x\in\Lambda_L}$ be a collection of vectors from $\E_\Omega$ satisfying that $\bw_x\ne0$ for only a finite number of~$x\in\Z^d$ and $\sum_{x\in\Lambda_L}\bw_x=0$. Let $\langle-\rangle_{J,\bb}^{(L)}$ denote the expectation with respect to $\nu_{J,\bb}^{(L)}$. Then we have the infrared bound~\cite{FSS,FILS1,FILS2},
\begin{equation}
\sum_{x,y\in\Lambda_L}\bigl\langle (\bw_x,\bS_x)(\bw_y,\bS_y)\bigr\rangle_{J,\bb}^{(L)}
\le J^{-1}\sum_{x,y\in\Lambda_L}(\bw_x,\bw_y)
\,D^{-1}_L(x,y)
\end{equation}
where
\begin{equation}
D^{-1}_L(x,y)=\frac1{|\Lambda_L^\star|}
\sum_{k\in\Lambda_L^\star\smallsetminus\{0\}}\frac1{\widehat D(k)}e^{ik(x-y)}.
\end{equation}
Now, let $\hate_1,\dots,\hate_n$ be an orthogonal basis in $\E_\Omega$ and choose $\bw_x=w_x\hate_\ell$, where $(w_x)_{x\in\Z^d}$ is such that 
$w_x\ne0$ only for a finite number of~$x\in\Z^d$ and
\begin{equation}
\label{w-sum}
\sum_{x\in\Z^d}w_x=0.
\end{equation} 
Passing to the limit $L\to\infty$ in such a way that $\nu_{J,\bb}^{(L)}$ converges to the state $\nu_{J,\bb}$, and then summing over $\ell=1,\dots,n$ gets us the bound
\begin{equation}
\label{IRB-w}
\sum_{x,y\in\Z^d}w_xw_y\,\bigl\langle(\bS_x,
\bS_y)\bigr\rangle_{J,\bb}\le nJ^{-1}\sum_{x,y\in\Z^d}w_xw_y \,D^{-1}(x,y).
\end{equation}
So far we have \eqref{IRB-w} only for $(w_x)$ with a finite support. But, using that fact that both quantities $D^{-1}(x,y)$ and $\langle(\bS_x,\bS_y)\rangle_{J,\bb}$ are uniformly bounded, \eqref{IRB-w} is easily extended to all absolutely-summable $(w_x)_{x\in\Z^d}$ (i.e., those satisfying $\sum_{x\in\Z^d}|w_x|<\infty$) which obey the constraint \eqref{w-sum}.

Let $(v_x)$ be as specified in the statement of the Lemma and let $a=\sum_{x\in\Z^d}v_x$. Fix $K$, let $\Lambda_K$ be as above and define $(w_x^{(K)})$ by
\begin{equation}
w_x^{(K)}=v_x-\frac a{|\Lambda_K|}\1_{\{x\in\Lambda_K\}}.
\end{equation}
Clearly, these $(w_x^{(K)})$ obey the constraint \eqref{w-sum}. Our goal is to recover \eqref{IRBrel} from \eqref{IRB-w} in the $K\to\infty$ limit. Indeed, plugging this particular $(w_x^{(K)})$ into \eqref{IRB-w}, the left hand side opens into four terms. The first of these is the sum of $v_xv_y\langle (\bS_x,\bS_y)\rangle_{J,\bb}$, which is part of what we want in \eqref{IRBrel}. The second and the third terms are of the same form and both amount to
\begin{equation}
a\sum_{x,y}v_x\1_{\{x\in\Lambda_K\}}\bigl\langle (\bS_x,\bS_y)\bigr\rangle_{J,\bb}
=a\Bigl\langle \sum_xv_x\Bigl(\bS_x,\frac1{|\Lambda_K|}\sum_{y\in\Lambda_K}\bS_y\Bigr)\Bigr\rangle_{J,\bb}.
\end{equation}
By our assumption of a sharp block-average magnetization in $\nu_{J,\bb}$, the average of the spins in $\Lambda_K$ can be replaced, in the $K\to\infty$ limit, by~$\bm$. Similarly, we claim that
\begin{equation}
\lim_{K\to\infty}\frac1{|\Lambda_K|^2}\sum_{x,y\in\Lambda_K}\bigl\langle (\bS_x,\bS_y)\bigr\rangle_{J,\bb}=|\bm|^2,
\end{equation}
so, recalling the definition of~$a$, the left hand side is in a good shape.

As for the right-hand side of \eqref{IRB-w} with $(w_x)=(w_x^{(K)})$, here we invoke the fact that (for $d\ge3$)
\begin{equation}
\lim_{K\to\infty}\frac1{|\Lambda_K|}\sum_{x\in\Lambda_L}D^{-1}(x,y)=0,
\end{equation}
uniformly in $y\in\Z^d$. The claim therefore follows.
\end{proofsect}

Next we show that for any parameters~$J$ and~$\bb$, and any $\bm_\star\in\mathscr{M}_\star(J,\bb)$, we can always find a state with magnetization $\bm_\star$ that is a limit of states satisfying Properties~1 and~2.

\begin{lemma}
\label{lemma2.1} For all\/ $J>0$, all\/ $\bb\in\E_\Omega$ and all\/\/ $\bm_\star\in\mathscr{M}_\star(J,\bb)$, there are sequences $(J_k)$, $(\bb_k)$ and $(\bm_k)$ with $J_k\to J$, $\bb_k\to\bb$, $\bm_k\to\bm_\star$ and $\mathscr{M}(J_k,\bb_k)=\{\bm_k\}$. In particular, there is a sequence $(\nu_{J_k,\bb_k})$ of infinite-volume Gibbs measures satisfying Properties~1 and~2, which weakly converge (possibly along a subsequence) to a measure $\nu_{J,\bb}$ with magnetization $\bm_\star$.
\end{lemma}

\begin{proofsect}{Proof}
The proof uses a little more of the convexity theory, let us recapitulate the necessary background. Let $f\colon\R^n\to(-\infty,\infty)$ be a convex and continuous function. Let $(\cdot,\cdot)$ denote the inner product in $\R^n$. For each~$x\in\R^n$, let $S(x)$ be the set of all possible limits of the gradients $\nabla f(x_k)$ for sequences~$x_k\in\R^n$ such that~$x_k\to x$ as $k\to\infty$. Then Theorem~25.6 of~\cite{Rocky} says that the set of all subgradients $\partial f(x)$ of $f$ at~$x$,
\begin{equation}
\partial f(x)=\bigl\{a\in\R^n\colon f(y)-f(x)\ge(y-x,a),\, y\in\R^n\bigr\},
\end{equation}
can be written as
\begin{equation}
\label{subdiff}
\partial f(x)=\conv(S(x))
\end{equation}
where $\conv(S(x))$ is the closed, convex hull of $S(x)$. (Here we noted that since the domain of~$f$ is all of~$\R^n$, the so called normal cone is empty at all~$x\in\R^n$.) But $S(x)$ is closed and thus $\conv(S(x))$ is simply the convex hull of $S(x)$. Now, by Corollary~18.3.1 of~\cite{Rocky}, we also know that if $S\subset\R^n$ is a bounded set of points and $C$ is its convex hull (no closure), then every extreme point of $C$ is a point from~$S$. Thus, we conclude: \emph{every extreme point of $\partial f(x)$ lies in~$S(x)$}.

Now we can apply the above general facts to our situation. Let $F(J,\bb)$ be the infinite-volume free energy of the model in \eqref{Ham}. Noting that $F(J,\bb)$ is defined for all $J\in\R$ and all $\bb\in\E_\Omega$, the domain of~$F$ is $\R\times\E_\Omega$. By well known arguments, $F$ is continuous and concave. Moreover, a comparison of \eqref{physF} and \eqref{subdiff} shows that $\mathscr{K}_\star(J,\bb)$ is---up to a sign change---the subdifferential of $F$ at $(J,\bb)$. As a consequence of the previous paragraph, every extreme point $[e_\star,\bm_\star]\in\mathscr{K}_\star(J,\bb)$ is given by a limit $\lim_{k\to\infty}[e_k,\bm_k]$, where $[e_k,\bm_k]$ are such that $\mathscr{K}_\star(J_k,\bb_k)=\{[e_k,\bm_k]\}$ for some $J_k\to J$ and $\bb_k\to\bb$. But $\bm_\star\in\mathscr{M}_\star(J,\bb)$ implies that $[e_\star,\bm_\star]$ is an extreme point of $\mathscr{K}_\star(J,\bb)$ for some $e_\star$, so the first part of the claim follows. 

To prove the second part, note that any infinite-volume limit of the finite-volume Gibbs state with periodic boundary condition and parameters $J_k$ and $\bb_k$ must necessarily have energy density $e_k$ and magnetization $\bm_k$. By compactness of the set of all Gibbs states (which is ensured by compactness of~$\Omega$), there is at least one (subsequential) limit $\langle-\rangle_{J,\bb}$ of the torus states as $J_k\to J$ and $\bb_k\to\bb$, which is then a translation-invariant Gibbs state with parameters~$J$ and~$\bb$ such that 
\begin{equation}
e_\star=\bigl\langle(\bS_x,\bS_y)\bigr\rangle_{J,\bb}
\quad\text{and}\quad\bm_\star=\langle\bS_x\rangle_{J,\bb},
\end{equation}
where~$x$ and $y$ is any pair of nearest neighbors of $\Z^d$. However, the block-average values of both quantities must be constant almost-surely, because otherwise $\langle-\rangle_{J,\bb}$ could have been decomposed into at least two ergodic states with distinct values of energy-density/magnetization pair, which would in turn contradict that
$[e_\star,\bm_\star]$ is an extreme point of $\mathscr{K}_\star(J,\bb)$.
\end{proofsect}

We note that the limiting measure is automatically $\Z^d$-translation and rotation invariant and, in addition, satisfies the block-average property. But, in the cases that are of specific interest to the present work (i.e., when $\mathscr{M}_\star(J,\bb)$ contains several elements), there is little hope that such a state is a torus state. Nevertheless, we can prove:

\begin{corollary}
\label{cor3.3}
Let $J\ge0$ and $\bb\in\E_\Omega$. Then for any $\bm_\star\in\mathscr{M}_\star(J,\bb)$, there exists a state~$\nu_{J,\bb}$ with (block-average) magnetization~$\bm_\star$ for which the infrared bound \eqref{IRBrel} holds. Moreover, the state~$\nu_{J,\bb}$ is $\Z^d$-translation and rotation invariant.
\end{corollary}

\begin{proofsect}{Proof}
For $J=0$ we obviously have a unique Gibbs state and the claim trivially holds. Otherwise, all of this follows from the weak convergence of the $\nu_{J_k,\bb_k}$ discussed above.
\end{proofsect}

\subsection{Proof of Main Theorem}
Now we have all the ingredients ready to prove Lemma~\ref{lemma1.3}:

\begin{proofsect}{Proof of Lemma~\ref{lemma1.3}}
Fix $\bm_\star\in\mathscr{M}_\star(J,\bb)$ and let~$\nu_{J,\bb}$ be the state described in Corollary~\ref{cor3.3}. To prove our claim, it just remains to choose $(v_x)$ as follows:
\begin{equation}
v_x=\begin{cases}
\frac1{2d},\qquad&\text{if }|x|=1,
\\
0,\qquad&\text{otherwise},
\end{cases}
\end{equation}
and recall the definition of $I_d$ from \eqref{Id}.
\end{proofsect}

Having established Lemma~\ref{lemma1.3}, we are ready to give the proof of the Key Estimate:

\begin{proofsect}{Proof of Key Estimate}
Let $J\ge0$ and $\bb\in\E_\Omega$. Let $\bm_\star\in\mathscr{M}_\star(J,\bb)$ and let $\langle-\rangle_{J,\bb}$ be the state satisfying \eqref{magav} and \eqref{magvar}. Our goal is to prove the bound \eqref{cor-bd}. To that end, let $\bm_0=\bm_0(\bS)$ denote the spatially averaged magnetization of the neighbors of the origin. The rotation symmetry of the state $\langle-\rangle_{J,\bb}$ then implies
\begin{equation}
\bigl\langle(\bS_x,\bS_0)\bigr\rangle_{J,\bb}=
\bigl\langle(\bm_0,\bS_0)\bigr\rangle_{J,\bb}.
\end{equation}
Next, conditioning on the spin configuration in the neighborhood of the origin, we use the DLR condition for the state $\langle-\rangle_{J,\bb}$ which results in
\begin{equation}
\bigl\langle(\bm_0,\bS_0)\bigr\rangle_{J,\bb}=
\bigl\langle(\bm_0,\nabla G(J\bm_0+\bb))\bigr\rangle_{J,\bb}.
\end{equation}
Finally, a simple calculation, which uses the fact that $\bm_\star=\langle\bS_0\rangle_{J,\bb}=\langle\bm_0\rangle_{J,\bb}=\langle\nabla G(J\bm_0+\bb)\rangle_{J,\bb}$, allows us to conclude that
\begin{multline}
\label{mmbd}
\qquad
\bigl\langle(\bm_0,\nabla G(J\bm_0+\bb))\bigr\rangle_{J,\bb}-|\bm_\star|^2\\
=\Bigl\langle\bigl(\bm_0-\bm_\star,
\nabla G(J\bm_0+\bb)-\nabla G(J\bm_\star+\bb)\bigr)\Bigr\rangle_{J,\bb}.
\qquad
\end{multline}

To proceed with our estimates, we need to understand the structure of the double gradient of function~$G(\bh)$. Recall the notation $\langle-\rangle_{0,\bh}$ for the single-spin state tilted by the external field $\bh$. Explicitly, for each measurable function~$f$ on~$\Omega$, we have $\langle f(\bS)\rangle_{0,\bh}=e^{-G(\bh)}\langle f(\bS)e^{(h,\bS)}\rangle_0$. Then the components of the double gradient correspond to the components of the covariance matrix of the vector-valued random variable $\bS$. In formal vector notation, for any $\ba\in\E_\Omega$, 
\begin{equation}
(\ba,\nabla)^2 G(\bh)=\bigl\langle (\ba,\bS-\langle\bS\rangle_{0,\bh})^2\bigr\rangle_{0,\bh}.
\end{equation}
Pick $\bh_0,\bh_1\in\E_\Omega$. Then we can write
\begin{equation}
\bigl(\bh_1-\bh_0,\nabla G(\bh_1)-\nabla G(\bh_0)\bigr)=\int_0^1\textd\lambda\,\, 
\Bigl\langle\bigl(\bh_1-\bh_0,\bS-\langle\bS\rangle_{0,\bh_\lambda}\bigr)^2
\Bigr\rangle_{0,\bh_\lambda},
\end{equation}
where $\bh_\lambda=(1-\lambda)\bh_0+\lambda\bh_1$. But the inner product on the right-hand side can be bounded using the Cauchy-Schwarz inequality, and since
\begin{equation}
\bigl\langle|\bS-\langle\bS\rangle_{0,\bh_\lambda}|^2\bigr\rangle_{0,\bh_\lambda}\le\max_{\bS\in\Omega}(\bS,\bS)=\kappa,
\end{equation}
we easily derive that
\begin{equation}
\bigl(\bh_1-\bh_0,\nabla G(\bh_1)-\nabla G(\bh_0)\bigr)\le\kappa|\bh_1-\bh_0|^2.
\end{equation}
This estimate in turn shows that the right-hand side of \eqref{mmbd} can be bounded by $\kappa J\langle|\bm_0-\bm_\star|^2\rangle_{J,\bb}$. But for this we have the bound from Lemma~\ref{lemma1.3}: $\langle|\bm_0-\bm_\star|^2\rangle_{J,\bb}\le nJ^{-1}I_d$. Putting all the previous arguments together, \eqref{cor-bd} follows.
\end{proofsect}

\begin{proofsect}{Proof of Main Theorem}
This now follows directly by plugging \eqref{cor-bd} into \eqref{genbd}.
\end{proofsect}

\section{Proofs of results for specific models}
\label{sec4}\medskip\noindent
By and large, this section is devoted to the specifics of the three models described in Section~\ref{sec2}. Throughout the entire section, we will assume that $\bb=\bzero$ and henceforth omit~$\bb$ from the notation. We begin with some elementary observations which will be needed in all three cases of interest but which are also of some general applicability. 

\vbox{
\subsection{General considerations}
\label{sec4.1}\nopagebreak\vspace{-6mm}\nopagebreak
\subsubsection{Uniform closeness to global minima}
We start by showing that, for the systems under study, the magnetization is \emph{uniformly} close to a mean-field magnetization. Let $\mathscr{M}_\MF(J)$ denote the set of all local minima of~$\varPhi_J$. Obviously, if we know that the actual magnetization comes close to minimizing the mean-field free energy, in must be close to a minimum or a ``near-minimum'' of this function. A useful measure of this closeness is the following: For $J\in[0,\infty]$ and $\vartheta>0$, we~let
\begin{equation}
\label{DJ}
D_J(\vartheta)=\sup
\Bigl\{\dist\bigl(\bm,\mathscr{M}_\MF(J)\bigr)\,\Big|\,
\bm\in\conv(\Omega),\,
\varPhi_J(\bm)<F_\MF(J)+\vartheta\Bigr\},
\end{equation}
where $F_\MF(J)$ denotes the absolute minimum of~$\varPhi_J$. However, to control the ``closeness'' we will have to make some assumptions about the behavior of the (local) minima of~$\varPhi_J$. An important property ensuring the desired uniformity in all three models under study is as follows:
}

\begin{Uniformity Property}
If $J\ge0$ and if $\bm\in\conv(\Omega)$ is a global minimum of~$\varPhi_{J}$, then there is an $\epsilon>0$ and a continuous function $\bm^\sharp\colon[J-\epsilon,J+\epsilon]\to\conv(\Omega)$ such that $\lim_{J'\to J}\bm^\sharp(J')=\bm$ and~$\bm^\sharp(J')$ is a local minimum of~$\varPhi_{J'}$ for all $J'\in[J-\epsilon,J+\epsilon]$.
\end{Uniformity Property}

In simple terms, the Uniformity Property states that every global minimum can be extended into a one-parameter family of local minima. Based on the Uniformity Property, we can state a lemma concerning the limit of $D_J(\vartheta)$ as $\vartheta\downarrow0$:

\begin{lemma}
\label{lemma4.1}
Suppose that~$\varPhi_J$ satisfies the above Uniformity Property. Then for all $J_0>0$,
\begin{equation}
\lim_{\vartheta\downarrow0}\sup_{0\le J\le J_0}D_J(\vartheta)=0.
\end{equation}
\end{lemma}

\begin{proofsect}{Proof}
This is essentially an undergraduate exercise in compactness. Indeed, if the above fails, then for some $\epsilon>0$, we could produce a sequence $\vartheta_k\downarrow0$ and $J_k\in[0,J_0]$ such that
\begin{equation}
D_{J_k}(\vartheta_k)\ge6\epsilon.
\end{equation}
This, in turn, implies the existence of $\bm_k\in\conv(\Omega)$ such that
\begin{equation}
\dist\bigl(\bm_k,\mathscr{M}_\MF(J_k)\bigr)\ge 3\epsilon
\quad\text{while}\quad
\varPhi_{J_k}(\bm_k)<F_\MF(J_k)+\vartheta_k.
\end{equation}
Let us use~$J$ and~$\bm$ to denote the (subsequential) limits of the above sequences. Using the continuity of~$\varPhi_J(\bm)$, to the right of the \emph{while} we would have~$\varPhi_J(\bm)=F_\MF(J)$ and~$\bm$ is thus a global minimum of~$\varPhi_J$. By our hypothesis, for each $k$ sufficiently large, there is a local minimum $\bm^\sharp(J_k)$ of~$\varPhi_{J_k}$ with $\bm^\sharp(J_k)$ converging to $\bm$ as $k\to\infty$. Since $\bm_k$ is also converging to $\bm$, the sequences $\bm_k$ and $\bm^\sharp(J_k)$ will eventually be arbitrary close. But that contradicts the bound to the left of the \emph{while}.
\end{proofsect}

\subsubsection{Monotonicity of mean-field magnetization}
For spin systems with an internal symmetry (which, arguably, receive an inordinate share of attention), the magnetization usually serves as an order parameter. In the context of mean-field theory, what would typically be observed is an interval $[0,J_\MF]$ where $\bm=\bzero$ is  the  global minimizer of~$\varPhi_J$, while for $J>J_\MF$, the function~$\varPhi_J$ is  minimized by a non-zero~$\bm$. This is the case for all three models under consideration. (It turns out that whenever $\langle\bS\rangle_0=0$, the unique global minimum of~$\varPhi_J$ for~$J$ sufficiently small is $\bm=\bzero$.) 

In order to prove the existence of a symmetry-breaking transition, we need to prove that the models under considerations have a unique point where the local minimum $\bm=\bzero$ ceases the status of a global minimum. This amounts to showing that, once the minimizer of~$\varPhi_J$ has been different from zero, it will never jump back to $\bm=\bzero$. In the mean-field theory with interaction~\eqref{Ham}, this can be proved using the monotonicity of the energy density; an analogous argument can be used to achieve the same goal for the corresponding systems on~$\Z^d$.

\begin{lemma}
\label{lemma4.1a}
Let $J_1<J_2$ and let $\bm_1$ be a global minimizer of~$\varPhi_{J_1}$ and $\bm_2$ a global minimizer of~$\varPhi_{J_2}$. Then $|\bm_1|\le|\bm_2|$. Moreover, if $J\mapsto\bm(J)$ is a differentiable trajectory of local minima,~then
\begin{equation}
\label{fider}
\frac\textd{\textd J} \,\varPhi_J
\bigl(\bm(J)\bigr)=-\frac12\bigl|\bm(J)\bigr|^2.
\end{equation}
\end{lemma}

\begin{proofsect}{Proof}
The identity \eqref{fider} is a simple consequence of the fact that, if~$\bm$ is a local minimum of~$\varPhi_J$, then~$\nabla\varPhi_J(\bm)=0$. To prove the first part of the claim, let $J,J'\ge0$ and let $\bm$ be a minimizer of~$\varPhi_J$. Let~$F_\MF(J)$ be the mean-field free energy. First we claim that
\begin{equation}
\label{FID}
F_\MF(J)-F_\MF(J')\ge -\frac{J-J'}2|\bm|^2.
\end{equation}
Indeed, since $F_\MF(J)=\varPhi_J(\bm)$, we have from the definition of~$\varPhi_J$ that
\begin{equation}
F_\MF(J)=-\frac{J-J'}2|\bm|^2+\varPhi_{J'}(\bm).
\end{equation}
Then the above follows using that~$\varPhi_{J'}(\bm)\ge F_\MF(J')$. Let $J_1<J_2$ and $\bm_1$ and $\bm_2$ be as stated. Then \eqref{FID} for the choice $J=J_2$, $J'=J_1$ and $\bm=\bm_2$ gives
\begin{equation}
\frac{F_\MF(J_2)-F_\MF(J_1)}{J_2-J_1}\ge-\frac12|\bm_2|^2.
\end{equation}
while \eqref{FID} for the choice $J=J_1$, $J'=J_2$ and $\bm=\bm_1$ gives
\begin{equation}
\frac{F_\MF(J_1)-F_\MF(J_2)}{J_1-J_2}\le-\frac12|\bm_1|^2.
\end{equation}
Combining these two bounds, we have $|\bm_1|\le|\bm_2|$ as stated.
\end{proofsect}

\subsubsection{One-component mean-field problems}
Often enough, the presence of symmetry brings along a convenient property that the multicomponent mean-field equation \eqref{MFeq} can be reduced to a one-component problem. Since this holds for all cases under consideration and we certainly intend to use this fact, let us spend a few minutes formalizing the situation. 

Suppose that there is a non-zero vector~$\bomega\in\E_\Omega$ such that~$\nabla G(h\bomega)$ is colinear with~$\bomega$ (and not-identically zero) for all~$h$. As it turns out, then also $\nabla S(m\bomega)$ is colinear with~$\bomega$, provided $m\bomega\in\conv(\Omega)$. Under these conditions, let us restrict both~$\bh$ and~$\bm$ to scalar multiples of~$\bomega$ and introduce the functions
\begin{equation}
\label{gs1D}
g(h)=|\bomega|^{-2}G(h\bomega)\quad\text{and}\quad s(m)=|\bomega|^{-2}S(m\bomega).
\end{equation}
The normalization by $|\bomega|^{-2}$ ensures that~$s(m)$ is given by the Legendre transform of $g(h)$ via the formula \eqref{S}. Moreover, the mean-field free-energy function~$\varPhi_J(m\bomega)$ equals the $|\bomega|^2$-multiple of the function
\begin{equation}
\label{phi1D}
\phi_J(m)=-\frac12Jm^2-s(m).
\end{equation}
The mean-field equation \eqref{MFeq} in turn reads
\begin{equation}
\label{MFE1D}
m=g'(Jm).
\end{equation}
In this one-dimensional setting, we can easily decide about whether a solution to \eqref{MFE1D} is a local minimum of $\phi_J$ or not just by looking at the stability of the solutions under iterations of \eqref{MFE1D}:

\begin{lemma}
\label{lemma4.2}
Let~$m$ be a solution to \eqref{MFE1D} and suppose~$\phi_J$ is twice continuously differentiable in a neighborhood of~$m$. If
\begin{equation}
\label{gderbd}
Jg''(Jm)<1
\end{equation}
then $m$ is a local minimum of $\phi_J$. Informally, only ``dynamically stable'' solutions to the (on-axis) mean-field equation can be local minima of $\phi_J$.
\end{lemma}

We remark that the term ``dynamically stable'' stems from the attempt to find solutions to \eqref{MFE1D} by running the iterative scheme $m_{k+1}=g'(Jm_k)$.

\begin{proofsect}{Proof}
Let~$h$ and~$m$ be such that $g'(h)=m$, which is equivalent to $h=s'(m)$. An easy calculation then shows that $g''(h)=-(s''(m))^{-1}$. Suppose now that $m$ is a solution to \eqref{MFE1D} such that \eqref{gderbd} holds. Then $h=Jm$ and from \eqref{gderbd} we have
\begin{equation}
s''(m)=-\bigl(g''(Jm)\bigr)^{-1}<-J.
\end{equation}
But that implies
\begin{equation}
\phi_J''(m)=-J-s''(m)>-J+J=0,
\end{equation}
and, using the second derivative test, we conclude that~$m$ is a local minimum of $\phi_J$.
\end{proofsect}

With Lemmas~\ref{lemma4.1}, \ref{lemma4.1a} and~\ref{lemma4.2} established, our account of the general properties is concluded and we can start discussing particular~models. What follows in the next three subsections are the three respective models laid out in order of increasing difficulty. Our repeated---and not particularly elegant---strategy will be to pound at the various models using internal symmetry as the mallet. The upshot is inevitably that at most one component becomes dominant while all other components act, among themselves, like a system at high temperature. Thus all subdominant components are equivalent and the full problem has been reduced to an effective scalar model.
In short, there are some parallels between the various treatments. However, somewhat to our disappointment, we have not been able to find a unified derivation covering ``all models of this~sort.''

\subsection{Potts model}
In order to prove Theorem~\ref{Pottsres}, we need to establish (rigorously) a few detailed properties of the mean-field free-energy function~\eqref{MFPotts}. In the view of \eqref{bmPotts} we will interchangeably use the notations $\bm$ and $(x_1,\dots,x_q)$ to denote the same value of the magnetization.

\begin{lemma}
\label{lemma3.1} 
Consider the~$q$-state Potts model with $q\ge3$. Let~$\varPhi_J$ be the mean-field free-energy function as defined in~\eqref{MFPotts}. If $\bm\in\conv(\Omega)$ is a local minimum of~$\varPhi_J$ then the corresponding $(x_1,\dots,x_q)$ is a permutation of the probability vector $(x_1^\star,\dots,x_q^\star)$ such that
\begin{equation}
x_1^\star\ge x_2^\star=\dots=x_q^\star.
\end{equation}
Moreover, when~$x_1^\star>x_2^\star$, we also have
\begin{equation}
\label{4.16b}
Jx_1^\star>1>Jx_2^\star.
\end{equation}
\end{lemma}

A complete proof of the claims in Lemma~\ref{lemma3.1} was, to our best knowledge, first provided in~\cite{KS}. (Strictly speaking, in~\cite{KS} it was only shown that the \emph{global} minima of~$\varPhi_J$ take the above form; however, the proof in~\cite{KS} can be adapted to also accommodate \emph{local} minima.) We will present a nearly identical proof but with a different interpretation of the various steps. The advantage of our reinterpretation is that it is easily applied to the other models of interest in this paper.

\begin{proofsect}{Proof of Lemma~\ref{lemma3.1}}
If~$\bm$ corresponds to the  vector $(x_1,\dots,x_q)$, we let~$\varPhi_J^{(q)}(x_1,\dots,x_q)$ be the quantity~$\varPhi_J(\bm)$. Suppose that $(x_1,\dots,x_q)$ is a local minimum. It is easy to verify that $(x_1,\dots,x_q)$ cannot lie on the boundary of $\conv(\Omega)$, so~$x_k>0$ for all $k=1,\dots,q$. Pick any two coordinates---for simplicity we assume that our choice is~$x_1$ and~$x_2$---and let $y=1-(x_3+\dots+x_q)$, $z_1=x_1/y$ and $z_2=x_2/y$. (Note that $y=x_1+x_2$ and, in particular, $y>0$.) Then~we~have
\begin{equation}
\varPhi_J^{(q)}(x_1,\dots,x_q)=-\frac12 Jy^2(z_1^2+z_2^2)+y(z_1\log
z_1+z_2\log z_2)+R_J^{(q)}(x_3,\dots,x_q),
\end{equation}
where $R_J^{(q)}(x_3,\dots,x_q)$ is independent of $z_1$ and $z_2$. Examining the form of the free energy, we find that the first two terms are proportional to the mean-field free-energy function of the Ising ($q=2$) system with reduced coupling $Jy$:
\begin{equation}
\label{Potts-Ising}
\varPhi_J^{(q)}(x_1,\dots,x_q)=y\,\varPhi_{Jy}^{(2)}(z_1,z_2)
+R_J^{(q)}(x_3,\dots,x_q).
\end{equation}

Since the only $z$-dependence is in the first term, the pair $(z_1,z_2)$ must be a local minimum of~$\varPhi_{Jy}^{(2)}$ regardless of what~$x_3,\dots,x_q$ look like. But this reduces the problem to the Ising model, about which much is known and yet more can easily be derived. The properties of~$\varPhi_J^{(2)}(z_1,z_2)$ we will need are:
\begin{enumerate}
\item[(i)]
$J_{\text{c}}=2$ is the critical coupling. For $J\le
J_{\text{c}}$, the free-energy function~$\varPhi_J^{(2)}(z_1,z_2)$ is lowest when $z_1=z_2$, while for $J>J_{\text{c}}$, the free-energy function~$\varPhi_J^{(2)}(z_1,z_2)$ is lowest when $\rho=|z_1-z_2|$ is the maximal (non-negative) solution to $\rho=\tanh(\frac12J\rho)$.
\item[(ii)]
Whenever $J>J_{\text{c}}$, the maximal solution to $\rho=\tanh(\frac12J\rho)$ satisfies $J(1-\rho^2)<2$, which implies that either $Jz_1>1$ and $Jz_2<1$ or \emph{vice versa}.
\item[(iii)]
For all~$J$ and $z_1\ge z_2$, the mean-field free-energy function~$\varPhi_J^{(2)}(z_1,z_2)$ monotonically decreases as $\rho=z_1-z_2$ moves towards the non-negative global minimum.
\end{enumerate}
All three claims are straightforward to derive, except perhaps~(ii), which is established by noting that, whenever $\rho>0$ satisfies the (Ising) mean-field equation, we have
\begin{equation}
\frac12J(1-\rho^2)=\frac
J{2\cosh(\frac12J\rho)^2}=\frac{J\rho}{\sinh(J\rho)}<1.
\end{equation}
Hence, if $J>J_{\text{c}}$ and $z_1> z_2$, then $Jz_2=\frac12J(1-\rho)<\frac12J(1-\rho^2)<1$ and thus $Jz_1>1$ because $J(z_1+z_2)=J>J_{\text{c}}=2$.

Based on (i-iii), we can draw the following conclusions for any pair of distinct indices~$x_j$ and~$x_k$: If $J(x_j+x_k)\le2$, then~$x_j=x_k$, because the $(k,j)$-th Ising pair is subcritical,
while if $J(x_j+x_k)>2$ then, using our observation~(ii), either $Jx_k>1$ and $Jx_j<1$ or \emph{vice versa}. But then we cannot have $Jx_k>1$ for more than one index $k$, because if
$Jx_k>1$ and $Jx_j>1$, we would have $J(x_j+x_k)>2$ and the $(k,j)$-th Ising pair would not be at a local minimum. All the other indices must then be equal because the associated two-component Ising systems are subcritical. Consequently, only one index from $(x_1,\dots,x_q)$ can take a larger value; the other indices are equal.
\end{proofsect}

\begin{proposition}
\label{prop4.2}
Consider the~$q$-state Potts model with $q\ge3$. Let~$\varPhi_J$ be the mean-field free-energy function as defined in~\eqref{MFPotts}. There there exist $J_1$ and $J_2=q$ with $J_1<J_2$ such that
\begin{enumerate}
\item[(1)]
$\bm=\bzero$ is a local minimum of~$\varPhi_J$ provided $J<J_2$.
\item[(2)]
$\bm=x_1^\star\hatv_1+\dots+x_q^\star\hatv_1$ with~$x_1^\star>x_2^\star=\dots=x_q^\star$ is a local minimum of~$\varPhi_J$ provided that $J>J_1$ and~$x_1^\star=\frac1q+m$, where $m$ is the maximal positive solution to the equation \eqref{mJPotts}.
\item[(3)]
For all $J\ge0$, there are no local minima except as specified in (1) and (2).
\end{enumerate}
Moreover, if $J_\MF$ is as in \eqref{JMF}, then the unique global minimum of~$\varPhi_J$ is as in~(1) for $J< J_\MF$ while for $J>J_\MF$ the function~$\varPhi_J$ has~$q$~distinct global minimizers as described in~(2) .
\end{proposition}

\begin{proofsect}{Proof of Proposition~\ref{prop4.2}}
Again, most of the above stated was proved in \cite{KS} but without the leeway for  local minima. (Of course, the formulas \eqref{mJPotts} and \eqref{JMF} date to an earlier epoch, see e.g.~\cite{Wu}.) What is not either easily derivable or already proved in \cite{KS} amounts to showing that if $m$ is a ``dynamically stable'' solution to \eqref{mJPotts}, the corresponding $\bm=x_1^\star\hatv_1+\dots+x_q^\star\hatv_1$ as described in~(2) is a local minimum for  the full~$\varPhi_J(\bm)$. The rest of this proof is spent proving the latter claim.

We first observe that for the set
\begin{equation}
\UU(x)=\bigl\{\bm=(x,x_2,\dots,x_q)\colon Jx_k\le1,\,k=2,\dots,q\bigr\}
\end{equation}
the unique (strict) global minimum of~$\varPhi_J$ occurs at 
\begin{equation}
\label{fmm}
\bm(x)=\bigl(x,\tfrac{1-x}{q-1},\dots,\tfrac{1-x}{q-1}\bigr).
\end{equation}
Indeed, otherwise we could further lower the value of~$\varPhi_J$ by bringing one of the $(j,k)$-th Ising pairs closer to its equilibrium, using the properties~(ii-iii) above. Now, suppose that $m$ satisfying \eqref{mJPotts} is ``dynamically  stable'' in the sense of Lemma~\ref{lemma4.2}. By \eqref{4.16b} we have that the corresponding~$x_1^\star=\frac1q+m$ satisfies $Jx_1^\star>1$ while the common value of~$x_k^\star$ for $k=2,\dots,q$ is such that $Jx_k^\star<1$. Suppose that the corresponding $\bm$ is \emph{not} a local minimum of the full~$\varPhi_J$. Then there exists a sequence~$(\bm_k)$ tending to~$\bm$ such that~$\varPhi_J(\bm_k)<\varPhi_J(\bm)$. But then there is also a sequence~$\bm_k'$ such that~$\varPhi_J(\bm_k')<\varPhi_J(\bm)$ where each~$\bm_k'$ now takes the form \eqref{fmm}. This contradicts that the restriction of~$\varPhi_J$ to the ``diagonal,'' namely the function $\phi_J(m)$, has a local minimum at $m$.
\end{proofsect}

Now we are ready to prove our main result about the~$q$-state Potts model.

\begin{proofsect}{Proof of Theorem~\ref{Pottsres}}
By well known facts from the FK representation of the Potts model, the quantities $e_\star(J)$ and $m_\star(J)$ arise from the pair $[e_\star^{\text{w}},\bm_\star^{\text{w}}]$ corresponding to the state with constant boundary conditions (the \emph{wired} state). Therefore, $[e_\star^{\text{w}},\bm_\star^{\text{w}}]$ is an extreme point of the convex set $\mathscr{K}_\star(J)$ and $\bm_\star^{\text{w}}\in\mathscr{M}_\star(J)$ for all~$J$. In particular, the bound \eqref{5} for $\bm_\star^{\text{w}}$ can be used without~apology.

Let $\delta_d$ be the part of the error bound in \eqref{5} which does not depend on~$J$. Explicitly, we have $\delta_d=\frac1{2q}(q-1)^2I_d$, because $\kappa=(q-1)/q$ and $\dim\,\E_\Omega=q-1$. Since $I_d\to0$ as $d\to\infty$, we have $\delta_d\to0$ as $d\to\infty$. Let us define
\begin{equation}
\label{e2def}
\epsilon_1=\epsilon_1(d,J)=\sup_{0\le J'\le J}D_{J'}(J\delta_d),
\end{equation}
where $D_J$ is as in \eqref{DJ}. It is easy to check that the Uniformity Property holds. Lemma~\ref{lemma4.1} then guarantees that every (extremal) physical magnetization $\bm_\star\in\mathscr{M}_\star(J)$ has to lie within $\epsilon_1$ from a local minimum~$\varPhi_J$. Since the asymmetric minima exist only for $J> J_1>0$ while $\bm=\bzero$ is a local minimum only for $J<J_2=q$, we have $m_\star(J)\le\epsilon_1$ for $J\le J_1$, while $|m_\star(J)-m_\MF(J)|\le\epsilon_1$ for $J>J_2$. But from the FKG properties of the random cluster representation we know that $J\mapsto m_\star(J)$ is non-decreasing so there must be a point, $J_{\text{\rm t}}\in(J_1,J_2]$, such that \twoeqref{Potts-mag-lb}{Potts-mag-ub} hold.

It remains to show that $|J_{\text{\rm t}}-J_\MF|$ tends to zero as $d\to\infty$. For $J\in(J_1,J_2)$, let $\varphi_{\text{S}}(J)$, resp., $\varphi_{\text{A}}(J)$ denote the value of~$\varPhi_J$ at the symmetric, resp., asymmetric local minima. The magnetization corresponding to the asymmetric local minimum exceeds some $\varkappa>0$ throughout $(J_1,J_2)$. Integrating \eqref{fider} with respect to~$J$ and using that $\varphi_{\text{S}}(J_\MF)=\varphi_{\text{A}}(J_\MF)$ then gives us~the~bound
\begin{equation}
\label{aaa}
\bigl|\varphi_{\text{S}}(J)-\varphi_{\text{A}}(J)\bigr|\ge\frac12\varkappa^2|J-J_\MF|.
\end{equation}
However, in the $\epsilon_1$-neighborhood $\UU_{\text{S}}(\epsilon_1)$ of the symmetric minimum, we will have
\begin{equation}
\label{aab}
\bigl|\varPhi_J(\bm)-\varphi_{\text{S}}(J)\bigr|\le\epsilon_1K,
\end{equation}
where $K$ is a uniform bound on the derivative of~$\varPhi_J(\bm)$ for $\bm\in\UU_{\text{S}}(\epsilon_1)$ and $J\in(J_1,J_2)$. Since the asymmetric minima are well separated from the boundary of $\conv(\Omega)$ for $J\in(J_1,J_2)$, a similar bound holds for the $\epsilon_1$-neighborhood of the asymmetric minimum. Comparing \twoeqref{aaa}{aab} and \eqref{5}, we find that if
\begin{equation}
\frac12\varkappa^2|J-J_\MF|-2\epsilon_1K>J\delta_d,
\end{equation}
no value of magnetization in the $\epsilon_1$-neighborhood of the local minima with a larger value of~$\varPhi_J$ is allowed. In particular, $|J_{\text{\rm t}}-J_\MF|\le \epsilon_2$ where $\epsilon_2=\epsilon_2(d)$ tends to zero as $d\to\infty$.
\end{proofsect}

\subsection{Cubic model}
Our first goal is to prove Proposition~\ref{prop-cubic}. We will begin by showing that the local minima of~$\varPhi_J$ and~$K^{(r)}_J$ are in one-to-one correspondence. Let us introduce the notation
\begin{equation}
X=\Bigl\{(\bar y,\bar\mu)\colon|\mu_j|\le1,\,y_j\ge0,\,\sum_{j=1}^ry_j=1\Bigr\}
\end{equation}
and let $X(\bm)$ denote the subspace of $X$ where $\bm=y_1\mu_1+\dots+y_r\mu_r$.

\begin{lemma}
\label{lemma4.3a}
Let $\bm\in\conv(\Omega)$ be a local minimum of~$\varPhi_J$. Then there exists a $(\bar y,\bar\mu)\in X(\bm)$ which is a local minimum of $K^{(r)}_J$ (as defined in~\eqref{MFcubic}).
\end{lemma}

\begin{proofsect}{Proof}
Let $\bm$ be a local minimum of~$\varPhi_J$. Since $X(\bm)$ is compact and $K^{(r)}_J$ is continuous on~$X$, the infimum
\begin{equation}
\label{4.27aa}
\varPhi_J(\bm)=\inf_{(\bar y,\bar\mu)\in X(\bm)}K^{(r)}_J(\bar y,\bar\mu)
\end{equation}
is attained at some $(\bar y,\bar\mu)\in X(\bm)$. We claim that this $(\bar y,\bar\mu)$ is a local minimum of $K^{(r)}_J$. Indeed, if the opposite is true, there is a sequence $(\bar y_k,\bar\mu_k)\in X$ converging to $(\bar y,\bar\mu)$ such that
\begin{equation}
\label{K1}
K^{(r)}_J(\bar y_k,\bar\mu_k)<K^{(r)}_J(\bar y,\bar\mu)=\varPhi_J(\bm).
\end{equation}
Now, $(\bar y,\bar\mu)$ was an absolute minimum of $K^{(r)}_J$ on $X(\bm)$, so $(\bar y_k,\bar\mu_k)\not\in X(\bm)$ and the magnetization~$\bm_k$ corresponding to $(\bar y_k,\bar\mu_k)$ is different from~$\bm$ for all~$k$. Noting that
\begin{equation}
\label{K2}
\varPhi_J(\bm_k)\le K^{(r)}_J(\bar y_k,\bar\mu_k)
\end{equation}
and combining \twoeqref{K1}{K2}, we thus have~$\varPhi_J(\bm_k)<\varPhi_J(\bm)$ for all~$k$. But~$\bm_k$ tends to~$\bm$ in $\conv(\Omega)$, which contradicts the fact that~$\bm$ is a local minimum of~$\varPhi_J$.
\end{proofsect}

Lemma~\ref{lemma4.3a} allows us to analyze the local minima in a bigger, simpler space:

\begin{lemma}
\label{lemma-cubic}
Let $K^{(r)}_J(\bar y,\bar\mu)$ be the quantity in \eqref{MFcubic}. Then each local minimum of $K^{(r)}_J(\bar y,\bar\mu)$ is an index-permutation of a state $(\bar y,\bar\mu)$ with $y_1\ge
y_2=\dots=y_r$ and $\mu_2=\dots=\mu_r=0$. Moreover, if $y_1>y_2$, then $\mu_1\ne0$.
\end{lemma}

\begin{proofsect}{Proof}
Let $(\bar y,\bar\mu)$ be a local minimum of
$K^{(r)}_J$ such that $y_1\ge y_2\ge\dots\ge y_r$ and fix a $k$ between~$1$ and~$r$. We abbreviate $y=y_k+y_{k+1}$ and introduce the variables $z_1=y_k/y$, $z_2=y_{k+1}/y$, $\nu_1=\mu_k$ and $\nu_2=\mu_{k+1}$. Then 
\begin{equation}
\label{kr2}
K^{(r)}_J(\bar y,\bar\mu)=y\,K^{(2)}_{Jy}(\bar z,\bar\nu)+R,
\end{equation}
where $K^{(2)}_{Jy}(\bar z,\bar\nu)$ is the mean-field free energy of an $r=2$ cubic model with coupling constant~$Jy$, and~$R$ is a quantity independent of $(\bar z,\bar\nu)$. As was mentioned previously, the $r=2$ cubic model is equivalent to two decoupled Ising models. Thus, 
\begin{equation}
\label{k2theta}
K^{(2)}_{Jy}(\bar z,\bar\nu)
=\varTheta_{Jy}(\rho_1)+\varTheta_{Jy}(\rho_2),
\end{equation}
where $\rho_1$ and $\rho_2$ are related to $z_1$, $z_2$, $\nu_1$ and $\nu_2$ via the equations
\begin{equation}
\label{z-nu}
\begin{array}{rclrcl}
z_1&=&{\textstyle\frac12}(1+\rho_1\rho_2), \qquad
z_1\nu_1&=&{\textstyle\frac12}(\rho_1+\rho_2),\\*[1mm]
z_2&=&{\textstyle\frac12}(1-\rho_1\rho_2), \qquad
z_2\nu_2&=&{\textstyle\frac12}(\rho_1-\rho_2).
\end{array}
\end{equation}
Now, the local minima of $\varTheta_J(\rho)$ occur at $\rho=\pm \rho(J)$, where $\rho(J)$ is the largest non-negative solution to the equation $\rho=\tanh(\frac12J\rho)$. Moreover, by the properties (i-iii) from the proof of Lemma~\ref{lemma3.1} we know that $\rho(J)=0$ for $J\le2$ while $\frac12J(1-\rho(J)^2)<1$ once $J>2$. From these observations we learn that if $y_k=y_{k+1}$, then $Jy\le2$ and $\mu_k=\mu_{k+1}=0$. On the other hand, if $y_k>y_{k+1}$, then $Jy>2$, $y_k=\frac12y(1+\rho(Jy)^2)$ and $y_{k+1}=\frac12y(1-\rho(Jy)^2)$ so, in particular, $Jy_k>1>Jy_{k+1}$. However, that forces that $k=1$, because otherwise we would also have $Jy_{k-1}>1$ and $J(y_{k-1}+y_k)>2$, implying that $(\bar y,\bar\mu)$ is not a local minimum of $K^{(r)}_J$ in the $(k-1,k)$-th sector. Hence, $y_2=\dots=y_r$ and $\mu_2=\dots=\mu_r=0$, while if $y_1>y_2$, then $\mu_1=\pm \rho(J)/z_1\ne0$.
\end{proofsect}

The proof of Lemma~\ref{lemma-cubic} gives us the following useful observation:

\begin{corollary}
\label{cor4.8}
Let $\bm=(m_1,m_2,\dots,m_r)$ be contained in $\conv(\Omega)$ and suppose that $m_1,m_2\ne0$. Then one of the four vectors
\begin{equation}
(m_1\pm m_2,0,m_3,\dots,m_r),\qquad (0,m_2\pm m_1,m_3,\dots,m_r)
\end{equation}
corresponds  to a magnetization $\bm'\in\conv(\Omega)$  with~$\varPhi_J(\bm')<\varPhi_J(\bm)$.
\end{corollary}

\begin{proofsect}{Proof}
Since $\bm$ is in the interior of $\conv(\Omega)$, there exists $(\bar y,\bar\mu)$ where  the infimum \eqref{4.27aa} is achieved. Let $z_1$, $z_2$, $\nu_1$ and $\nu_2$ be related to $y_1$, $y_2$, $\mu_1$ and $\mu_2$ as in \twoeqref{kr2}{z-nu}. Now by \eqref{k2theta} the free  energy  of the corresponding sector of $(\bar y,\bar\mu)$ equals the sum of the free energies of two decoupled Ising models  with biases $\rho_1$ and $\rho_2$. Without loss of generality, suppose that $\rho_1>\rho_2\ge0$. Recalling the property~(iii) from the proof of Lemma~\ref{lemma3.1}, $\rho\mapsto\varTheta_J(\rho)$ decreases when $\rho\ge0$ gets closer to the non-negative local minimum. Thus, if  $\rho_1$ is nearer to the local minimum of $\varTheta_{Jy}$ than $\rho_2$, by increasing $\rho_2$ we lower the free energy by a non-trivial amount. Similarly, if $\rho_2$ is the one that is closer, we decrease $\rho_1$. 

By inspection of \eqref{z-nu}, the former operation produces a new quadruple $z_1'$, $z_2'$, $\nu_1'$ and $\nu_2'$, with $\nu_2'=0$ and $z_1'\nu_1'=\rho_1$. But that corresponds to the magnetization vector $(m_1',m_2',m_3,\dots,m_r)$, where
\begin{equation}
m_1'=\rho_1 y=m_1+m_2\quad\text{and}\quad m_2'=0,
\end{equation}
which is what we stated above. The other situations are handled analogously.
\end{proofsect}

Now we are finally ready to establish the claim about local/global minima of~$\varPhi_J$:

\begin{proofsect}{Proof of Proposition~\ref{prop-cubic}}
By Lemma~\ref{lemma4.3a}, every local minimum of~$\varPhi_J$ corresponds to a local minimum of $K^{(r)}_J$. Thus, using Lemma~\ref{lemma-cubic} we know that all local minima $\bm$ of~$\varPhi_J$ will have at most one non-zero component. Writing $\bomega=(1,0,\dots,0)$, $\bh=h\bomega$ and $\bm=m\bomega$, we can use the formalism from Section~\ref{sec4.1}. In particular, the on-axis moment generating function $g(h)$ is~given~by
\begin{equation}
g(h)=-\log(2r)+\log(r-1+\cosh h).
\end{equation}
Differentiating this expression, \eqref{MFE1D} shows that every local minimum $m$ has to satisfy the equation~\eqref{cubic-eq}. Now, for $r>2$, a little work shows that $h\mapsto g'(h)$ is convex for
\begin{equation}
(r-1)^2-(r-1)\cosh h+2>0
\end{equation}
and concave otherwise. In particular, for $r>3$, the equation \eqref{cubic-eq} has either one non-negative solution $m=0$ or three non-negative solutions, $m=0$, $m=m_-(J)$ and $m=m_+(J)$, where $0\le m_-(J)\le m_+(J)$. However, $m_+(J)$ is ``dynamically stable'' and, using Lemma~\ref{lemma4.2}, $m_-(J)$ never corresponds to a local minimum. 

To finish the proof we need to show that $\bm=(m_+(J),0,\dots,0)$ is a local minimum of the full~$\varPhi_J$. If the contrary were true, we would have a sequence $\bm_k$ tending to $\bm$ such that~$\varPhi_J(\bm_k)<\varPhi_J(\bm)$. Then an ($r-1$)-fold use of Corollary~\ref{cor4.8} combined with the symmetry of~$\varPhi_J$ implies the existence of a sequence $\bm_k'=(m_k,0,\dots,0)$ tending to~$\bm$ and satisfying~$\varPhi_J(\bm_k')\le\varPhi_J(\bm_k)$ for all~$k$. But that contradicts that~$m_+(J)$ is a local minimum of the on-axis mean-field free energy function. So~$\bm$ was a local minimum of~$\varPhi_J$ after all. The existence of a unique mean-field transition point $J_\MF$ is a consequence of Lemma~\ref{lemma4.1a} and the fact that $\bm=\bzero$ ceases to be a local minimum for $J\ge r$.
\end{proofsect}

\begin{proofsect}{Proof of Theorem~\ref{Cubicres}}
The proof is basically identical to that of Theorem~\ref{Pottsres}, so we will be rather sketchy. First we note that $m_\star(J)$ is achieved at some extremal translation-invariant state whose magnetization $\bm_\star$ is an element of $\mathscr{M}_\star(J)$. Let $\delta_d=\frac12rI_d$  and define $\epsilon_1$ as in \eqref{e2def}. Then $\bm_\star$ has to be within $\epsilon_1$ from a local minimum of~$\varPhi_J$. While this time we cannot proclaim that $J\mapsto m_\star(J)$ is non-decreasing, all the benefits of monotonicity can be achieved by using the monotonicity of the energy density $e_\star(J)$. Indeed, $J\mapsto e_\star(J)$ is non-decreasing and, by Corollary~\ref{cor-ener1} and the Key Estimate, we have
\begin{equation}
\label{EMrel}
\Bigl|\,e_\star(J)-\frac12 m_\star(J)^2\Bigr|\le \frac J2rI_d=J\delta_d.
\end{equation}
But then $e_\star(J)$ must undergo a unique large jump at some $J_{\text{\rm t}}$ from values $e_\star(J)\le 2J\delta_d$ to values near $\frac12m_\MF(J)^2$ by less than $2J\delta_d$. So~$m_\star(J)$ has to jump at $J=J_{\text{\rm t}}$ as well, in order to obey~\eqref{EMrel}. The width of the ``transition region'' is controlled exactly as in the case of the Potts model.
\end{proofsect}

\subsection{Nematic model}
The nematic models present us with the difficulty that an explicit formula for~$\varPhi_J(\bm)$ seems impossible to derive. However, the situation improves in the dual Legendre variables. Indeed,
examining \twoeqref{S}{Fi}, it is seen that the stationary points of~$\varPhi_J(\bm)$ are in one-to-one correspondence with the stationary points of the (Gibbs) free-energy function
\begin{equation}
\varPsi_J(\bh)=\frac1{2J}|\bh|^2-G(\bh),
\end{equation}
via the relation $\bh=J\bm$. (In the case at hand, $\bh$ takes values in~$\E_\Omega$ which was defined as the space of all $N\times N$ traceless matrices.) Moreover, if $\bm=\nabla G(\bh)$, then we have
\begin{equation}
\varPsi_J(\bh)-\varPhi_J(\bm)=\frac1{2J}|\bh-J\bm|^2
\end{equation}
so the values $\varPsi_J(\bm)$ and~$\varPhi_J(\bh)$ at the corresponding stationary points are the same. Furthermore, some juggling with Legendre transforms shows that if~$\bm$ is a local minimum of~$\varPhi_J$, then $\bh=J\bm$ is a local minimum of~$\varPsi_J$. Similarly for local maxima and saddle points of~$\varPhi_J$. 

\begin{lemma}
\label{lemma4.7}
Each stationary point of\/ $\varPsi_J(\bh)$ on $\E_\Omega$
is a traceless $N\times N$ matrix~$\bh$ with eigenvalues that can be reordered to the form $h_1\ge h_2=\dots=h_N$.
\end{lemma}

\begin{proofsect}{Proof}
The claim is trivial for $N=2$ so let $N\ge3$. Without loss of generality, we can restrict ourselves to diagonal, traceless matrices~$\bh$. Let~$\bh=\text{diag}(h_1,\dots,h_N)$ be such that $\sum_\alpha h_\alpha=0$ and let $v_\alpha$, with $\alpha=1,\dots,N$, be the components a unit vector in~$\R^N$. Let $\langle-\rangle_0$ be the expectation with respect to the \emph{a priori} measure~$\mu$ on $\Omega$ and let $\langle-\rangle_\bh$ be the state on~$\Omega$ tilted by~$\bh$. Explicitly, we have
\begin{equation}
\label{4.40}
\langle f\rangle_\bh=e^{-G(\bh)}\int \mu(\textd\bv) f(\bv)\exp\Bigl\{\sum_{\alpha=1}^N h_\alpha v_\alpha^2\Bigr\}
\end{equation} 
for any measurable function $f$ on the unit sphere in~$\R^N$.

As in the case of the Potts and cubic models, the proof will be reduced to the two-component problem. Let $\bh$ be a stationary point of $\Psi_J$ and let $\alpha$ and $\beta$ be two distinct indices between $1$ and~$N$. The relevant properties of $\langle-\rangle_\bh$ are then as follows:
\begin{enumerate}
\item[(i)]
If $J\langle v_\alpha^4+v_\beta^4\rangle_\bh>3$, then $h_\alpha\ne h_\beta$.
\item[(ii)]
If $h_\alpha>h_\beta$, then $J\langle v_\alpha^4\rangle_\bh>\frac32>J\langle v_\beta^4\rangle_\bh$.
\end{enumerate}
The proof of these facts involves a non-trivial adventure with modified Bessel functions, $I_n(x)$, where~$n$ is any non-negative integer and $I_n(x)=\frac1\pi\int_0^\pi\textd\theta\, e^{x\cos\theta}\cos(n\theta)$. To keep the computations succinct, we introduce the polar coordinates, $v_\alpha = r\cos\theta$ and $v_\beta = r\sin\theta$, where $\theta\in[0,2\pi)$ and $r\ge0$. Let $\langle-\rangle_{\alpha,\beta}$ denote the expectation with respect to the $r$-marginal of the state $\langle-\rangle_{\bh'}$ where $\bh'=\diag(h_1',\dots,h_N')$ is related to $\bh$ via $h_\alpha'=h_\beta'=\frac12(h_\alpha+h_\beta)$, while $h_\gamma'=h_\gamma$ for $\gamma\ne\alpha,\beta$. Explicitly, if $\bar f(r,\theta)$ corresponds to $f(v_\alpha,v_\beta)$ via the above change of coordinates,~then
\begin{equation}
\bigl\langle f(v_\alpha,v_\beta)\bigr\rangle_\bh=\frac{
\Bigl\langle\int_0^{2\pi} \textd\theta\, e^{r^2\Delta\cos(2\theta)}\,\bar f(r,\theta)\Bigr\rangle_{\alpha\beta}
}{\Bigl\langle\int_0^{2\pi}\textd\theta
\, e^{r^2\Delta\cos(2\theta)}\Bigr\rangle_{\alpha\beta}},
\end{equation}
where $\Delta=\frac12(h_\alpha-h_\beta)$.

We begin by deriving several identities involving modified Bessel functions.
First, a straightforward calculation shows that
\begin{equation}
\langle v_\alpha^2-v_\beta^2\rangle_\bh=A_{\alpha\beta}(\Delta)\,\bigl\langle r^2 I_1(r^2\Delta)\bigr\rangle_{\alpha\beta},
\end{equation}
where $A_{\alpha\beta}(\Delta)^{-1}=\langle I_0(r^2\Delta)\rangle_{\alpha\beta}$.
Similarly we get
\begin{equation}
\label{4.19}
\langle
v_\alpha^2v_\beta^2\rangle_\bh=A_{\alpha\beta}(\Delta)\,
\bigl\langle \tfrac18r^4 \bigl(I_0(r^2\Delta)-I_2(r^2\Delta)\bigr)\bigr\rangle_{\alpha\beta}.
\end{equation}
But $I_0(x)-I_2(x)=(2/x)I_1(x)$, whereby we have the identity
\begin{equation}
\label{miracle} 2(h_\alpha-h_\beta)\langle
v_\alpha^2v_\beta^2\rangle_\bh =\langle
v_\alpha^2-v_\beta^2\rangle_\bh.
\end{equation}
A similar calculation using trigonometric formulas shows that
\begin{align}
\label{4.21a}
\langle v_\alpha^4\rangle_\bh &=A_{\alpha\beta}(\Delta)\,
\bigl\langle r^4 \bigl(\tfrac38I_0(r^2\Delta)+\tfrac12I_1(r^2\Delta)+\tfrac18
I_2(r^2\Delta)\bigr)\bigr\rangle_{\alpha\beta},
\\
\label{4.21b}
\langle v_\beta^4\rangle_\bh &=A_{\alpha\beta}(\Delta)\,
\bigl\langle r^4 \bigl(\tfrac38I_0(r^2\Delta)-\tfrac12I_1(r^2\Delta)+\tfrac18
I_2(r^2\Delta)\bigr)\bigr\rangle_{\alpha\beta}.
\end{align}
In particular, since $I_0(0)=1$ while $I_1(0)=I_2(0)=0$, we have
\begin{equation}
\label{zazrak}
h_\alpha=h_\beta\quad \Rightarrow\quad \langle v_\alpha^4\rangle_\bh=\langle v_\beta^4 \rangle_\bh=3\langle v_\alpha^2v_\beta^2 \rangle_\bh.
\end{equation}
The identities \twoeqref{4.19}{zazrak} will now allow us to prove~(i-ii).

First we note that he fact that $\bh$ was a stationary point of $\varPsi_J$ implies that 
$h_\gamma-h_{\gamma'} =J\langle v_\gamma^2-v_{\gamma'}^2\rangle_\bh$ for all $\gamma,\gamma'=1,\dots,N$.
Using this in \eqref{miracle}, we have the following dichotomy
\begin{equation}
\label{jeste}
\text{either}\quad h_\alpha=h_\beta\quad\text{or}\quad
2J\langle v_\alpha^2v_\beta^2 \rangle_\bh=1.
\end{equation}
To establish (i), suppose that $J\langle v_\alpha^4+v_\beta^4\rangle_\bh>3$ but $h_\alpha=h_\beta$. Then \eqref{zazrak} gives us $2J\langle v_\alpha^2v_\beta^2 \rangle_\bh>1$, in contradiction with \eqref{jeste}. Hence, (i) must hold. To prove (ii), assume that $h_\alpha>h_\beta$ and note that then $\Delta>0$. Applying that $I_1(x)>0$ and $I_2(x)>0$ for~$x>0$ in \eqref{4.21a}, we easily show using \eqref{4.21a} that $\langle v_\alpha^4\rangle_\bh>3\langle v_\alpha^2v_\beta^2\rangle_\bh$. Similarly, the bound $I_1(x)>I_2(x)$ for~$x>0$, applied in \eqref{4.21b}, shows that $\langle v_\beta^4\rangle_\bh<3\langle v_\alpha^2v_\beta^2\rangle_\bh$. From here (ii) follows by invoking \eqref{jeste}.

Now we are ready to prove the desired claim. Let $\bh$ be a stationary point. First let us prove that there are no three components of $\bh$ such that $h_\alpha>h_\beta>h_\gamma$. Indeed, if that would be the case, (i-ii) leads to a contradiction, because $h_\alpha>h_\beta$ would require that $J\langle v_\beta^4\rangle_\bh<3/2$ while $h_\beta>h_\gamma$ would stipulate that $J\langle v_\beta^4\rangle_\bh>3/2$! Thus, any stationary point $\bh$ of $\varPsi_J$ can only have two values for $\langle v_\alpha^4\rangle_\bh$. However, if (say) both $\langle v_1^4\rangle_\bh$ and $\langle v_2^4\rangle_\bh$ take on the larger value (implying that $h_1=h_2$), then $J\langle v_1^4+v_2^4\rangle_\bh>3$ and $\bh$ cannot be a stationary point. From here the claim follows.
\end{proofsect}

The symmetry of the problem at hand allows us to restrict ourselves to the on-axis formalism from Section~\ref{sec4.1}. In particular, we let $\bomega=\diag(1,-\frac1{N-1},\dots,-\frac1{N-1})$, $\bh=h\bomega$ and $\blambda=\lambda\bomega$ and define the functions $g(h)$, $s(\lambda)$ and $\phi_J(\lambda)$ as in \twoeqref{gs1D}{phi1D}. Lemma~\ref{lemma4.7} in turn guarantees that all local minimizers of~$\varPhi_J$ appear within the domain of $\phi_J$. What remains to be proved is the converse. This can be done using some of the items established above.

\begin{lemma}
Suppose that~$\lambda$ is a stationary point of the scalar free energy~$\phi_J$ 
which satisfies $Jg''(J\lambda)<1$. Then $\blambda=\lambda\bomega$, with $\bomega=\diag(1,-\frac1{N-1},\dots,-\frac1{N-1})$, is a local minimizer of~$\varPhi_J$.
\end{lemma}

\begin{proofsect}{Proof}
To simplify the exposition, we will exploit the $O(N)$-symmetry of the problem: If $\bg\in O(N,\R)$ is any $N\times N$ orthogonal matrix, then
\begin{equation}
\varPhi_J(\bm)=\varPhi_J(\bg^{-1}\bm\,\bg),
\end{equation}
with similar considerations applying to $\varPsi_J(\bh)$.
Thus, for all intents and purposes, we may assume that the arguments of these functions are already in the diagonal form and regard the diagonal as an $N$-component vector.
(Indeed, we will transfer back and forth between the vector and matrix language without further ado.)

Again we are forced to work with the dual variables.
To that end, let $\psi_J(h)$ be the quantity $|\bomega|^{-2}\varPsi_J(h\bomega)$. Clearly, the relation between~$\psi_J$ and~$\phi_J$ is as for~$\varPsi_J$ and~$\varPhi_J$. First, let us demonstrate that every stationary point of the scalar free energy~$\psi_J$ represents a stationary point of the full~$\varPsi_J$. Indeed, let~$\K$ be the orthogonal complement of vector $\bomega$ in $\R^N$. As a simple computation shows, any $\bk\in\K$ has a zero first component. If $\bk=(0,k_2,\dots,k_N)\in\K$ is small, then
\begin{equation}
G(h\bomega+\bk)=G(h\bomega)+\Bigl\langle\sum_\beta k_\beta\, v_\beta^2\Bigr\rangle_{h\bomega}+O\bigl(|\bk|^2\bigr),
\end{equation}
where $\langle-\rangle_\bh$ is as in \eqref{4.40}. Now $\langle v_\beta^2\rangle_{h\bomega}$ is the same for all $\beta=2,\dots,N$, and in the view of the fact that $\sum_\beta k_\beta=0$, the expectation vanishes. Hence, $\nabla\varPsi_J(h\bomega)$ has all components corresponding to the subspace~$\K$ equal to zero. Now if~$h$ is a stationary point of~$\psi_J$, we know that $(\bomega,\nabla\varPsi_J(h\bomega))=0$ and thus $\nabla\varPsi_J(h\bomega)=\bzero$ as claimed.

To prove the desired claim, it now suffices to show that the Hessian of $\varPsi_J$ is positive definite at~$\bh=h^\star\bomega$ when~$h^\star$ satisfies~$Jg''(h^\star)<1$. (Recall that the corresponding stationary points of~$\psi_J$ and~$\phi_J$ are related by $h=J\lambda$.) This in turn amounts to showing that $\nabla\nabla G(h\bomega)$ is dominated by the $J^{-1}$-multiple of the unit matrix. Although we must confine ourselves to~$\E_\Omega$, it is convenient to consider the Hessian of $G(\bh)$ in a larger space which contains the constant vector and restrict our directional probes to vectors from~$\E_\Omega$. In general, the entries of the Hessian are given in terms of truncated correlation functions:
\begin{equation}
\bigl(\text{Hess}(G)\bigr)_{\alpha\beta}=\langle v_\alpha^2v_\beta^2\rangle_\bh-\langle v_\alpha^2\rangle_\bh\langle v_\beta^2\rangle_\bh.
\end{equation}
For the problem at hand, there are only four distinct entries:
\newcommand{\sA}{{\text{\rm\small A}}}
\newcommand{\sB}{{\text{\rm\small B}}}
\newcommand{\sC}{{\text{\rm\small C}}}
\newcommand{\sD}{{\text{\rm\small D}}}
\begin{equation}
\text{Hess}(G)=\left(
\begin{array}{ccccc}
\sA & \sB & \dots & \dots & \sB \\
\sB & \sC & \sD & \dots & \sD \\
\vdots & \sD & \ddots & \ddots & \vdots \\
\vdots & \vdots & \ddots & \sC & \sD \\
\sB & \sD & \dots & \sD & \sC \\
\end{array}
\right).
\end{equation}
Clearly, $\bomega$ itself is an eigenvector of~$\text{Hess}(G)$ with the eigenvalue $\sA-\sB$. On the other hand, if $\bk\in\K$, then the first row and column of~$\text{Hess}(G)$ are irrelevant. Writing the remaining $(N-1)\times(N-1)$ block in the form $(\sC-\sD)\1+\sC\,\BbbS$, where~$\BbbS$ is the matrix with \emph{all} entries equal to one, it follows easily that all of~$\K$ is an eigenspace of~$\text{Hess}(G)$ with eigenvalue~$\sC-\sD$.

It remains to show that these eigenvalues are strictly smaller than~$J^{-1}$. The first one, namely, $\sA-\sB$ is less than $J^{-1}$ by our assumption that $Jg''(h^\star)<1$. As to the other eigenvalue, $\sC-\sD$, we note that
\begin{equation}
\sC-\sD = \langle v_\alpha^4\rangle_\bh-\langle v_\alpha^2v_\beta^2\rangle_\bh,\qquad \alpha>\beta>1.
\end{equation}
Now, equation \eqref{zazrak} tells us that, under our conditions, $\langle v_\alpha^2v_\beta^2\rangle_\bh$ \emph{equals} $\frac13\langle v_\alpha^4\rangle_\bh$. So we need that $\frac23\langle v_\alpha^4\rangle_\bh$ is less than~$J$. But since $h_1=h^\star>h_\alpha$, that is exactly the condition~(ii) derived in the proof of Lemma~\ref{lemma4.7}.
\end{proofsect}

Now we are ready to establish our claims concerning the local minima of~$\varPhi_J$:

\begin{proofsect}{Proof of Proposition~\ref{prop2.4a}}
Let $\bomega$ be as above and note that $|\bomega|^2=N/(N-1)$. Then the on-axis moment generating function from \eqref{gs1D} becomes
\begin{equation}
\label{gnematic}
g(h)=\frac{N-1}N\,\log\int\pi_N(\textd\bv)\,e^{h\frac N{N-1}(v_1^2-\frac1N)},
\end{equation}
where~$\pi_N$ is the uniform probability measure on the unit sphere in~$\R^N$ and~$v_1$ is the first component of~$\bv$. An argument involving the~$N$-dimensional spherical coordinates then shows that
\begin{equation}
\pi_N(v_1\in\textd x)=C(N)\,(1-x^2)^{\frac{N-3}2}\textd x,
\end{equation}
where $C(N)$ is the ratio of the surfaces of the unit spheres in~$\R^{N-1}$ and~$\R^N$. By substituting this into \eqref{gnematic} and applying \eqref{MFE1D}, we easily find that, in order for $\blambda=\lambda\bomega$ to be a local minimum of~$\varPhi_J$, the scalar $\lambda$ has to satisfy the equation \eqref{lMFeq}.

A simple analysis of \eqref{lMFeq} shows that for $J\ll1$, the only solution to \eqref{lMFeq} is $\lambda=0$, while for $J\gtrsim N^2$, the solution $\lambda=0$ is no longer perturbatively stable. Since Lemma~\ref{lemma4.1a} guarantees that the norm of all global minimizers increases with~$J$, there must be a unique $J_\MF\in(0,\infty)$ and a non-decreasing function $J\mapsto\lambda_\MF(J)$ such that $\lambda_\MF(J)$ solves \eqref{lMFeq} and that every global minimizer of~$\varPhi_J$ at any $J>J_\MF$ which is a continuity point of $J\mapsto\lambda_\MF(J)$ corresponds to $\lambda=\lambda_\MF(J)$. (At any possible point of discontinuity of $J\mapsto\lambda_\MF(J)$, the $\lambda$ corresponding to any global minimizer is sandwiched between $\lim_{J'\uparrow J}\lambda_\MF(J')$ and $\lim_{J'\uparrow J}\lambda_\MF(J')$.) The claim is thus~proved.
\end{proofsect}

In order to prove the large-$N$ part of our statements concerning the mean-field theory of the nematic model, we will need to establish the following scaling property:

\begin{lemma}
\label{lemma4.8}
Let~$\varPhi_J^{(N)}$ denote the free-energy function of the $O(N)$-nematic Hamiltonian. Introduce the matrix $\bomega=\diag(1,-\frac1{N-1},\dots,-\frac1{N-1})$ and define the normalized mean-field free-energy function
\begin{equation}
\phi^{(N)}_J(\lambda)=\frac1N |\bomega|^{-2}\varPhi_{JN}^{(N)}(\lambda\bomega),
\qquad \lambda<1.
\end{equation}
Then, as $N\to\infty$, the function $\lambda\mapsto\phi^{(N)}_J(\lambda)$ converges, along with all of its derivatives, to the function
\begin{equation}
\phi^{(\infty)}_J(\lambda)=-\frac J2\lambda^2+\frac12\log\frac1{1-\lambda}.
\end{equation}
\end{lemma}

\begin{proofsect}{Proof}
The proof is a straightforward application of Laplace's method to the measure on the right-hand side of \eqref{lMFeq}. Indeed, for any $h\ge0$, consider the measure~$\rho_{h,N}$ on~$[0,1]$ defined by
\begin{equation}
\rho_{h,N}(\textd x)=\frac{(1-x^2)^{\frac{N-3}2}e^{h Nx^2}}
{\int_0^1\textd x\,(1-x^2)^{\frac{N-3}2}e^{h Nx^2}}\,\textd x.
\end{equation}
Noting that the function~$x\mapsto (1-x^2)^{\frac12}e^{h x^2}$ has a unique maximum at~$x=x_h$, where
\begin{equation}
x_h^2=\max\Bigl\{0,1-\frac1{2h}\Bigr\},
\end{equation}
we easily conclude that
\begin{equation}
\lim_{N\to\infty}\rho_{h,N}(\cdot)=\delta_{x_h}(\cdot),
\end{equation}
where $\delta_a(\cdot)$ denotes the Dirac point mass at~$x=a$. Here the limit taken in the sense of weak convergence on the space of all bounded continuous functions on~$[0,1]$.
The proof of this amounts to standard estimates for the Laplace method; we leave the details to the reader.

Let $g_N(h)$ denote the function $g(hN)$ where $g$ is as in~\eqref{gnematic}. Since any derivative of $g_N(h)$ can be expressed as a truncated correlation function of measure $\rho_{h,N}$, we easily conclude that $h\mapsto g_N(h)$ converges, along with all of its derivatives, to the function
\begin{equation}
\label{ginf}
g_\infty(h)=\lim_{N\to\infty}g_N(h)=\max\Bigl\{0,h-\frac12-\frac12\log(2h)\Bigr\},
\end{equation}
for all $h\ge0$. Now, the function $s_N(\lambda)=\frac1N|\bomega|^{-2} S(\lambda\bomega)$---where $S(\cdot)$ is the entropy of the $O(N)$-nematic model---is the Legendre transform of $g_N$, so we also get
\begin{equation}
s_\infty(\lambda)=\lim_{N\to\infty}s_N(\lambda)=-\frac12\log\frac1{1-\lambda}.
\end{equation}
(Again, the convergence extends to all derivatives, provided $\lambda<1$.)
From here the claim follows by noting that $\phi_J^{(N)}(\lambda)=-\frac J2\lambda^2-s_N(\lambda)$, which tends to $\phi^{(\infty)}_J(\lambda)$ in the desired sense.
\end{proofsect}

\begin{proofsect}{Proof of Proposition~\ref{prop2.4b}}
By Lemma~\ref{lemma4.8}, the scaled mean-field free-energy function~$\phi_J^{(N)}$ is, along with any finite number of its derivatives, uniformly close to~$\phi_J^{(\infty)}$ on compact subsets of $[0,1)$, provided~$N$ is sufficiently large. Now the local minima of~$\phi_J^{(\infty)}$ will again satisfy a mean-field equation, this time involving the function~$g_\infty$ from \eqref{ginf}. Since \begin{equation}
g'(h)=\begin{cases}
1-\frac1{2h},\qquad&\text{if }h>\frac12,
\\
0,\qquad&\text{otherwise},
\end{cases}
\end{equation}
there are at most two perturbatively stable solutions to the mean-field equation: One at $\lambda=0$ and the other at
\begin{equation}
\lambda=\frac12\bigl(1+\sqrt{1-4J^{-2}}\bigr).
\end{equation}
Moreover, these local minima interchange the role of the global minimum at some finite and non-zero~$J^{(\infty)}_\MF$, which is a solution of a particular transcendental equation. For~$J$ near~$J^{(\infty)}_\MF$, the second derivative of~$\phi_J^{(\infty)}$ is uniformly positive around both local minima.

The convergence stated in Lemma~\ref{lemma4.8} ensures that all of the previously listed facts will be (at least qualitatively) satisfied by~$\phi_J^{(N)}$ for~$N$ large as well. Thus,~$\phi_J^{(N)}$ has at most one positive local minimum, which immediately implies that $J\mapsto\lambda_\MF^{(N)}(J)$ is continuous whenever it is defined. Moreover, since the local minima of~$\phi_J^{(N)}$ converge to those of~$\phi_J^{(\infty)}$, we also easily recover the asymptotic statements~\twoeqref{limlam}{limJ}. This finishes the proof.
\end{proofsect}

\begin{proofsect}{Proof of Theorem~\ref{nematicres}}
The proof is similar to that of the Potts and cubic models; the only extra impediment is that now we cannot take for granted that there is only one non-zero local minimum. As before, most of the difficulties will be resolved by invoking the monotonicity of the energy density~$e_\star(J)$, which is defined e.g.\ by optimizing $\frac12\langle(\bQ_0,\bQ_x)\rangle_J$ over all Gibbs states invariant under the lattice translations and rotations. 

In the present case, $\kappa$ and $n$ in the Main Theorem are given by $\kappa=(N-1)/N$ and $n=\frac12N(N-1)$. Thus, letting $\delta_d=\frac14(N-1)^2I_d$, the quantity $J\delta_d$ is the corresponding error term on the right-hand side of~\eqref{5}. Define $\epsilon_1$ by the formula \eqref{e2def}. Then Lemma~\ref{lemma4.7} guarantees that the diagonal form~$\blambda$ of~$\langle\bQ_0\rangle_J$ for any Gibbs state is an index permutation of a vector of the type
\begin{equation}
\label{vektor}
\Bigl(\lambda+a_1,-\frac\lambda{N-1}+a_2,\dots,-\frac\lambda{N-1}+a_N\Bigr),
\end{equation}
where $\sum_ia_i=0$, $\sum_ia_i^2\le\epsilon_1^2$ and $\lambda$ corresponds to a local minimum of~$\varPhi_J$. If $\blambda$ is the physical magnetization giving rise to $\lambda_\star(J)$, we let $\lambda_\MF^\star(J)$ be a value of $\lambda$, corresponding to a local minimum of~$\varPhi_J$, for which $\blambda$ takes the form \eqref{vektor}. Then Corollary~\ref{cor-ener1} and the Key Estimate~give
\begin{equation}
\Bigl|e_\star(J)-\frac12\frac N{N-1}\lambda_\MF^\star(J)^2\Bigr|\le 2J\delta_d.
\end{equation}
Now for $J\le J_0\ll1$, we know the only local minimum is for $\lambda_\MF^\star(J)=0$, while for $J\ge J_1\gtrsim N^2$, the zero vector is no longer a local minimum and hence $\lambda_\MF^\star(J)$ exceeds some $\varkappa'>0$. But $J\mapsto e_\star(J)$ is non-decreasing so there must be a $J_{\text{\rm t}}\in[J_0,J_1]$ where~$e_\star(J)$ jumps by at least $\varkappa'-2J_{\text{\rm t}}\delta_d$, which is positive once~$d$ is sufficiently large.
The fact that $J_{\text{\rm t}}$ must be close to $J_\MF$ for large enough~$d$ is proved exactly as for the Potts and cubic models.
\end{proofsect}

\section{Mean-field theory and complete-graph models}
\label{sec5}\medskip\noindent
Here we will show that the mean-field formalism developed in Section~\ref{sec1.2} has a very natural interpretation for the model on a complete graph. An important reason for the complete graph picture is to provide a tangible physical system to motivate some of the physical arguments. The forthcoming derivation is a rather standard exercise in large-deviation theory~\cite{Ellis,Dembo-Zeitouni}, so we will keep it rather brief.

We will begin by a precise definition of the problem. Let $\GG_N$ be a complete graph on~$N$ vertices and consider a spin system on $\GG_N$ with single-spin space $\Omega$ and the Hamiltonian
\begin{equation}
\beta H_N(\bS)=-\frac JN\sum_{1\le x<y\le N}(\bS_x,\bS_y)-\sum_{x=1}^N(\bb,\bS_x).
\end{equation}
(Recall that $\Omega$ is a compact subset of a finite-dimensional vector space $\E_\Omega$ with inner product denoted as in the previous formula.) Let $\mu$ denote the \emph{a priori} spin measure and let $\langle-\rangle_0$ denote the corresponding expectation. For each configuration $\bS$, introduce the empirical magnetization by the formula
\begin{equation}
\bm_N(\bS)=\frac1N\sum_{x=1}^N\bS_x.
\end{equation}
If $\bm\in\conv(\Omega)$ and $\epsilon>0$, let $\UU_\epsilon(\bm)$ denote the $\epsilon$-neighborhood of $\bm$ in $\conv(\Omega)$ in the metric induced by the inner product on $\E_\Omega$. Then we have:

\begin{theorem}
\label{thm5.1}
For each $\bm\in\conv(\Omega)$,
\begin{equation}
\label{LDP}
\lim_{\epsilon\downarrow0}\lim_{N\to\infty}\frac1N\log\Bigl\langle
e^{-\beta H_N(\bS)}\1_{\{\bm_N(\bS)\in\,\UU_\epsilon(\bm)\}}\Bigr\rangle_0=
-\varPhi_{J,\bb}(\bm),
\end{equation}
where~$\varPhi_{J,\bb}(\bm)$ is as defined in Section~\ref{sec1.2}. Moreover, if $\nu_N$ denotes the Gibbs measure obtained by normalizing $e^{-\beta H_N(\bS)}$ and if $F_\MF(J,\bb)$ denotes the infimum of~$\varPhi_{J,\bb}(\bm)$ over $\bm\in\conv(\Omega)$, then
\begin{equation}
\lim_{N\to\infty}\nu_N\bigl(\varPhi_{J,\bb}
(\bm_N(\bS))\ge F_\MF(J,\bb)+\epsilon\bigr)=0 
\end{equation}
for every $\epsilon>0$.
\end{theorem}

\begin{proofsect}{Proof}
By our assumption, $\E_\Omega$ is a finite-dimensional vector space. Moreover, $\Omega$ is compact and thus the logarithmic generating function $G(\bh)$ defined in \eqref{G} exists for all $\bh\in\E_\Omega$. As a consequence of Cram\'er's Theorem for i.i.d. random variables on $\R^n$, see Theorem~2.2.30 in~\cite{Dembo-Zeitouni}, the measures
\begin{equation}
\mu_N(\cdot)=\mu\bigl(\bm_N(\bS)\in\cdot\bigr)
\end{equation}
satisfy a large-deviation principle on $\R^d$ with rate function \eqref{S}. In particular,
\begin{equation}
\lim_{\epsilon\downarrow0}\lim_{N\to\infty}
\frac1N\log\mu_N\bigl(\UU_\epsilon(\bm)\bigr)
=S(\bm),\qquad \bm\in\conv(\Omega).
\end{equation}
Now $\beta H_N$ can be written as follows
\begin{equation}
\beta H_N=N E_{J,\bb}\bigl(\bm_N(\bS)\bigr)-\frac JN\sum_{x=1}(\bS_x,\bS_x).
\end{equation}
Since the second term is bounded by a non-random constant almost surely and since $\bm\mapsto E_{J,\bb}(\bm)$ is uniformly continuous throughout $\conv(\Omega)$, \eqref{LDP} follows by inspecting the definition of~$\varPhi_{J,\bb}(\bm)$.
\end{proofsect}

\section*{Acknowledgments}
\noindent
The research of L.C.~was supported by the NSF under the grant DMS-9971016
and by the NSA under the grant NSA-MDA~904-00-1-0050.

\end{document}

\vspace{1cm} \scshape\footnotesize

\hbox to \textwidth{ \hglue 0.6 cm \vbox{ \hbox{Marek Biskup}
\hbox{Department of Mathematics} \hbox{UCLA} \hbox{Los Angeles CA
90095-1555} \hbox{U.S.A.} \hbox{email: \rm biskup@math.ucla.edu} }
\hfill \vbox{ \hbox{Philippe Blanchard \& Tyll Kr\"uger}
\hbox{Department of Theoretical Physics} \hbox{University of
Bielefeld, BiBoS} \hbox{Postfach 100131, D-33501} \hbox{Bielefeld,
Germany} \hbox{email: \rm
$\{$blanchard,$\,$tkrueger$\}$@physik.uni-bielefeld.de} } \hglue
0.3 cm }
\end{document}